\documentclass{article}

\newif\ifarxiv
\arxivtrue

\usepackage{graphicx}
\usepackage{amsmath,amssymb,amsthm,bm}
\usepackage[a4paper]{geometry}
\usepackage{algorithm,algorithmic}
\usepackage[colorlinks=false,pdfborder={0 0 0}]{hyperref}
\usepackage{authblk}
\usepackage{color}
\usepackage{todonotes}
\usepackage{makecell}

\graphicspath{{fig/}}
\DeclareGraphicsExtensions{.mps,.pdf,.eps,.png}

\newenvironment{keywords}{\medskip\textbf{Keywords:}}{}
\newenvironment{MSCcodes}{\medskip\textbf{MSC codes:}}{}

\newtheorem{theorem}{Theorem}

\newtheorem{example}{Example}

\newtheorem{remark}{Remark}
\newtheorem{lemma}{Lemma}

\theoremstyle{plain}

\DeclareMathOperator*{\argmin}{\arg\min}

\DeclareMathOperator{\Span}{span}

\newcommand{\fro}{\mathsf F}

\newcommand*{\trans}{^{\top}}

\newcommand{\bmat}[1]{\begin{bmatrix}#1\end{bmatrix}}

\newcommand*{\abs}[1]{\bigl\lvert#1\bigr\rvert}
\newcommand*{\norm}[1]{\bigl\Vert#1\bigr\rVert}
\newcommand*{\normbig}[1]{\big\Vert#1\big\rVert}
\newcommand*{\normBig}[1]{\Big\Vert#1\Big\rVert}
\newcommand*{\normbigg}[1]{\bigg\Vert#1\bigg\rVert}
\newcommand*{\normF}[1]{\bigl\Vert#1\bigr\rVert_{\fro}}
\newcommand*{\normFbig}[1]{\big\Vert#1\big\rVert_{\fro}}
\newcommand*{\normFBig}[1]{\Big\Vert#1\Big\rVert_{\fro}}
\newcommand*{\normFbigg}[1]{\bigg\Vert#1\bigg\rVert_{\fro}}

\def\adots{\mathinner{\mkern2mu\raise1pt\hbox{.}\mkern2mu
    \raise4pt\hbox{.}\mkern2mu\raise7pt\hbox{.}\mkern1mu}}

\newcommand*{\macheps}{\bm u}

\newcommand*{\epsqr}{\omega_{\mathrm{qr}}}

\newcommand*{\bigO}{O}

\newcommand{\sigmin}{\sigma_{\min}}

\newcommand*{\epsAv}{\varepsilon_{AZ}}
\newcommand*{\epsB}{\delta_{M_R^{-1}}}

\newcommand*{\epsH}{\omega_{H}}
\newcommand*{\epskappaW}{\omega_{\kappa(W)}}
\newcommand*{\epskapparW}{\omega_{\kappa(r, W)}}
\newcommand*{\epsls}{\varepsilon_{ls}}
\newcommand*{\epsMLAZ}{\delta_{M_L^{-1}AZ}}
\newcommand*{\epsMLb}{\delta_{M_L^{-1}b}}

\newcommand*{\epsMV}{\varepsilon_{mv}}

\newcommand*{\epsr}{\varepsilon_{r}}

\newcommand*{\krylov}{\mathcal{K}}
\newcommand*{\ki}{i^*}
\newcommand*{\kj}{j^*}

\title{On the backward stability of s-step GMRES}
\author[1]{Erin Carson}
\author[1]{Yuxin Ma}
\affil[1]{Department of Numerical Mathematics, Faculty of Mathematics and Physics, Charles University, Sokolovsk\'{a} 49/83, 186 75 Praha 8, Czechia

{\tt Email: carson@karlin.mff.cuni.cz, yuxin.ma@matfyz.cuni.cz}}

\begin{document}
\maketitle

\begin{abstract}
Communication, i.e., data movement, is a critical bottleneck for the performance of classical Krylov subspace method solvers on modern computer architectures. Variants of these methods which avoid communication have been introduced, which, while equivalent in exact arithmetic, can be unstable in finite precision. 
In this work, we address the backward stability of \(s\)-step GMRES, also known as communication-avoiding GMRES.
Building upon the ``modular framework'' proposed in [A.~Buttari, N.~J.~Higham, T.~Mary, \& B.~Vieubl{\'e}. Preprint in 2024.], we present an improved framework for simplifying the analysis of \(s\)-step GMRES, which includes standard GMRES (\(s=1\)) as a special case, by isolating the effects of rounding errors in the QR factorization and the solution of the least squares problem.
The key advantage of this new framework is that it is evident how the orthogonalization method affects the backward error, and it is not necessary to re-evaluate anything other than the orthogonalization itself when modifying the orthogonalization used in GMRES.
Using this framework, we analyze \(s\)-step GMRES with popular block orthogonalization methods: block modified Gram--Schmidt and reorthogonalized block classical Gram--Schmidt algorithms.

An example illustrates the resulting instability of \(s\)-step GMRES when paired with the classical \(s\)-step Arnoldi process and shows the limitations of popular strategies for resolving this instability.
To address this issue, we propose a modified \(s\)-step Arnoldi process that allows for much larger block size \(s\) while maintaining satisfactory accuracy, as confirmed by our numerical experiments.

\begin{keywords}
backward stability, \(s\)-step GMRES, communication-avoiding, Arnoldi process
\end{keywords}

\begin{MSCcodes}
65F10, 65F50, 65G50
\end{MSCcodes}
\end{abstract}

\section{Introduction}
\label{sec:introduction}

Given a nonsingular square matrix \(A \in \mathbb{R}^{n \times n}\) and a right-hand side \(b \in \mathbb R^{n}\), this work considers the iterative solution of the linear system
\begin{equation} \label{problem:linear-system}
    A x = b.
\end{equation}
A popular method for this problem is the generalized minimal residual algorithm (GMRES) introduced in~\cite{SS1986}, which chooses \(x^{(i)} \in x^{(0)} + \krylov_i(A, r)\) to minimize \(\norm{A x^{(i)} - b}_2\) in the \(i\)-th iteration with a Krylov subspace \(\krylov_i(A, r) = \Span\{r, A r, \dotsc, A^{i-1} r\}\) and the initial residual \(r = b - A x^{(0)}\).

The standard GMRES algorithm creates and orthogonalizes the Krylov basis one vector at a time, predominantly making use of BLAS 2 operations.
On modern computer architectures, the performance of classical iterative solvers is heavily dominated by communication, i.e., data movement and synchronization.
This motivated the introduction of \(s\)-step (also called communication-avoiding) variants of the GMRES algorithm, which can reduce asymptotic latency costs by a factor of \(s\) and can take advantage of BLAS 3 operations. 
In~\cite[Section 3.6]{Hoemmen2010}, it is demonstrated that CA-GMRES provides a speedup ranging from \(1.31\times\) to \(4.3\times\) compared to the standard GMRES on an Intel platform with 8 cores. 
The authors in~\cite{YATHD2014} show that CA-GMRES can achieve speedups of up to \(2\times\) over the standard GMRES on one or multiple GPUs. Particularly, the orthogonalization phase in CA-GMRES can reach up to \(4.16\times\) speedup compared to the same phase in the standard GMRES when \(s = 10\) is used.

In each iteration, the \(s\)-step GMRES algorithm generates a block Krylov basis consisting of \(s\) vectors formed as follows:
\begin{equation}
    \bmat{p_0(A) v& p_1(A) v& \dotsc & p_{s-1}(A) v},
\end{equation}
where \(p_0\), \(\dotsc\), \(p_{s-1}\) are specified polynomials, and then performs a step of block orthogonalization.
For each iteration of the standard GMRES algorithm, \(s=1\) and \(p_0(A)=I\).
The \(s\)-step GMRES algorithm is equivalent in exact arithmetic to the standard GMRES algorithm, however, it is known that they can behave quite differently in finite precision for \(s > 1\). 
In~\cite{Hoemmen2010, IE2017}, the authors conjecture that the stability of \(s\)-step GMRES depends on the condition number of the block Krylov basis, which is affected by the choice of the polynomial.
Although monomials are a natural choice for the polynomial basis, \cite{Beckermann1997} demonstrated that the condition number of the block Krylov basis increases exponentially with block size \(s\) using a monomial basis, which has empirically been observed to have a negative effect on stability. 
Alternative polynomials, such as Newton and Chebyshev polynomials, have been suggested to mitigate the growth rate of the condition number with \(s\), as detailed in \cite{Carson2015, Hoemmen2010}.
Additionally, the authors in \cite{IE2017} propose adaptively varying \(s\) across different iterations of \(s\)-step GMRES also to mitigate the increase in the condition number of the basis.

Considering rounding errors, the backward stability of the standard GMRES algorithms with Householder and modified Gram--Schmidt orthogonalization (MGS) has been investigated in~\cite{DGRS1995, PRS2006}.
Recently, \cite{BHMV2024} introduced a comprehensive framework to simplify the rounding error analysis of GMRES algorithms.
Building on this framework, we analyze the rounding error of the \(s\)-step GMRES algorithm to illustrate how the condition number of the basis affects the backward error in this study.
Our analysis is presented with the goals of making it possible to easily identify the sources of errors, clarifying the impact of different orthogonalization methods on the backward error of the computed $s$-step GMRES results, and broadening the applicability of the framework introduced by~\cite{BHMV2024}.
In addition, our analysis formally shows that the backward error of \(s\)-step GMRES is influenced by the condition number of the basis.
From this analysis, we further propose a modified Arnoldi process to allow for the use of a substantially larger block size \(s\), as confirmed by numerical experiments, increasing the cost of the orthogonalization by a factor of 2.

The remainder of this paper is organized as follows. In Section \ref{sec:preliminary}, we give an overview of the $s$-step Arnoldi and GMRES algorithms. 
In Section~\ref{sec:stability},
we present our abstract framework for the analysis of the $s$-step GMRES method. 
In Section \ref{sec:ortho}, 
we invoke our abstract framework to analyze the backward error of the \(s\)-step GMRES algorithm with different commonly-used block orthogonalization methods. We discuss the implications of the theoretical results and comment on stopping criteria in Section \ref{sec:discussion}.
Then we propose a modified $s$-step Arnoldi process in Section~\ref{sec:modifiedArnoldi} to improve the backward stability of the \(s\)-step GMRES.
Numerical experiments are presented in Section~\ref{sec:experiments} which compare the \(s\)-step GMRES algorithm with the modified and classical $s$-step Arnoldi processes.

We first introduce some notation used throughout the paper. We use MATLAB indexing to denote submatrices. 
For example, we use \(X_{1:i}\) to denote the first \(i\) columns of 
\(X\), and use \(H_{1:i+1, 1:i}\) to denote \((i+1)\)-by-\(i\) leading submatrix of \(H\).
For simplicity, we also omit the column indices of the square submatrices.
For instance, we abbreviate \(i\times i\) submatrix of \(R\) as \(R_{1:i} := R_{1:i, 1:i}\).
In addition, we use superscripts to denote the iteration.
For instance, \(x^{(i)}\) denotes the approximate solution in the \(i\)-th iteration.
We use \(\norm{\cdot}\) to denote the $2$-norm and \(\normF{\cdot}\) to denote the Frobenius norm in bounds,
and we further use \(\kappa(A)\) to represent the 2-norm condition number defined by \(\norm{A}/\sigmin(A)\), where \(\sigmin(A)\) is the smallest singular value of \(A\).
We also use \(\hat{\cdot}\) to denote computed quantities, and \(\macheps\) to represent the unit roundoff.
For polynomials, we indicate the degree of polynomial using subscripts.

\section{The \(s\)-step Arnoldi and GMRES algorithms}
\label{sec:preliminary}
In this section, we begin with an overview of the \(s\)-step Arnoldi and GMRES algorithms without employing any preconditioners, followed by a discussion of their variants incorporating left and right preconditioners.

The standard GMRES algorithm, namely the \(s\)-step GMRES algorithm with \(s = 1\), utilizes the Arnoldi process to build an orthonormal basis \(V_{1:i+1}\) for the Krylov subspace \(\krylov_{i+1}(A, r)\), where \(r = b - A x^{(0)}\).
This process is expressed by
\begin{equation*}
    \bmat{r & W_{1:i}} = V_{1:i+1} R_{1:i+1}
\end{equation*}
where \(W_{1:i} = A V_{1:i}\) and \(R_{1:i+1}\) is an upper triangular matrix.
After \(i\) iterations, if convergence is attained, the solution is updated via \(x^{(i)} \gets x^{(0)} + V_{1:i} y^{(i)}\) with \(y^{(i)} = \argmin_y \norm{\beta e_1 - H_{1:i+1, 1:i} y}\).
Here, \(\beta = \norm{r}\), \(e_1\) represents the first column of the identity matrix of size \((i+1) \times (i+1)\), and \(H_{1:i+1, 1:i} = R_{1:i+1, 2:i+1}\).

When \(s > 1\), the \(s\)-step GMRES algorithm handles \(s\) vectors at a time in every iteration.
Specifically, in the \(i\)-th iteration, the \(s\)-step algorithm first builds an orthonormal basis \(V_{1:is+1}\), by the \(s\)-step Arnoldi process, satisfying
\begin{equation*}
    \bmat{r & W_{1:is}} = V_{1:is+1} R_{1:is+1},
\end{equation*}
where \(W_{1:is} = A B_{1:is}\) and each subblock of \(B_{1:is}\), namely \(B_{(k-1)s+1:ks}\) with \(k = 1, 2, \ldots, i\), satisfies
\begin{equation*}
    B_{(k-1)s+1:ks} = \bmat{p_0(A) V_{(k-1)s+1}& p_1(A) V_{(k-1)s+1}& \dotsb & p_{s-1}(A) V_{(k-1)s+1}}
\end{equation*}
with given polynomials \(p_0\), \(p_1\), \(\ldots\), \(p_{s-1}\).
Next, the \(s\)-step GMRES algorithm updates \(x^{(i)}\) by \(x^{(0)} + B_{1:is} y^{(i)}\), where \(y^{(i)} = \argmin_y\norm{\beta e_1 - H_{1:is+1, 1:is} y}\) with \(H_{1:is+1, 1:is} = R_{1:is+1, 2:is+1}\).
Notice that only the subdiagonal elements are nonzero in the lower triangular part of \(H_{1:is+1, 1:is}\).
The least squares problem \(\min_y\norm{\beta e_1 - H_{1:is+1, 1:is} y}\) is often solved by applying the Givens QR factorization to \(H_{1:is+1, 1:is}\) such that
\begin{equation*}
    H_{1:is+1, 1:is} = G_{1:is+1} T_{1:is+1, 1:is}
\end{equation*}
with an orthogonal matrix \(G_{1:is+1}\) and upper triangular matrix \(T_{1:is}\).
Then \(y^{(i)}\) can be computed by solving the triangular system \(T_{1:is} y^{(i)} = \beta G_{1, 1:is}\trans\), since 
\begin{align*}
    y^{(i)} &= \argmin_y\norm{\beta e_1 - G_{1:is+1} T_{1:is+1, 1:is} y} \\
    &= \argmin_y\norm{\beta G_{1, 1:is+1}\trans - T_{1:is+1, 1:is} y} \\
    &= \argmin_y\norm{\beta G_{1, 1:is}\trans - T_{1:is} y}.
\end{align*}
This \(s\)-step GMRES algorithm is the so-called ``non-traditional'' variant introduced by~\cite{IE2017}.

Furthermore, if we consider the left and right preconditioners, i.e., \(M_L\) and \(M_R\), \(s\)-step Arnoldi process aims to build an orthonormal basis \(V_{1:is+1}\) of \(\krylov_{i+1} (M_L^{-1} A M_R^{-1}, r)\) with \(r = M_L^{-1}(b - A x^{(0)})\).
This means that \(V_{1:is+1}\) satisfies
\begin{equation*}
    \bmat{r & W_{1:is}} = V_{1:is+1} R_{1:is+1}
\end{equation*}
with \(W_{1:is} = M_L^{-1} A M_R^{-1} B_{1:is}\).
Then \(s\)-step GMRES updates \(x^{(i)}\) through \(x^{(0)} + M_R^{-1} B_{1:is} y^{(i)}\), where \(y^{(i)}\) is the optimal solution for \(\min_y\norm{\beta e_1 - H_{1:is+1, 1:is} y}\) with \(\beta = \norm{r}\) and \(H_{1:is+1, 1:is} = R_{1:is+1, 2:is+1}\).

We summarize the above preconditioned \(s\)-step Arnoldi and GMRES algorithms in Algorithms~\ref{alg:sstep-Arnoldi-classical} and~\ref{alg:sstep-GMRES}, respectively.

\begin{algorithm}[!tb]
\begin{algorithmic}[1]
    \caption{The \(i\)-th iteration of the \(s\)-step Arnoldi process \label{alg:sstep-Arnoldi-classical}}
    \REQUIRE
    A matrix \(A \in \mathbb R^{n\times n}\), a vector \(r\), a block size \(s\), a left-preconditioner \(M_L \in \mathbb R^{n\times n}\), a right-preconditioner \(M_R \in \mathbb R^{n\times n}\),
    the basis \(B_{1:(i-1)s}\) and the preconditioned basis \(Z_{1:(i-1)s}\) generated by the first \(i-1\) classical Arnoldi steps,
    the matrix \(W_{1:(i-1)s}\), the orthonormal matrix \(V_{1:(i-1)s+1}\), and the upper triangular matrix \(R_{1:(i-1)s+1}\) satisfying \(\bmat{r & W_{1:(i-1)s}} = V_{1:(i-1)s+1} R_{1:(i-1)s+1}\).
    \ENSURE
    The basis \(B_{1:is}\), the preconditioned basis \(Z_{1:is}\), the matrices \(W_{1:is}\), \(V_{1:is+1}\), and \(R_{1:is+1}\) satisfying \(\bmat{r & W_{1:is}} = V_{1:is+1} R_{1:is+1}\).

    \STATE \(B_{(i-1)s+1:is} \gets \bmat{p_0(A) V_{(i-1)s+1} & p_1(A) V_{(i-1)s+1} & \dotsi & p_{s-1}(A) V_{(i-1)s+1}}\). \label{line:Arnoldi:B}
    \STATE \(Z_{(i-1)s+1:is} \gets M_R^{-1} B_{(i-1)s+1:is}\)
    \STATE \(W_{(i-1)s+1:is} \gets M_L^{-1} A Z_{(i-1)s+1:is}\).
    \STATE Compute the QR factorization of \(\bmat{r & W_{1:is}} = V_{1:is+1} R_{1:is+1}\) based on \(\bmat{r & W_{1:(i-1)s}} = V_{1:(i-1)s+1} R_{1:(i-1)s+1}\). \label{line:classical-orth}
\end{algorithmic}
\end{algorithm}

\begin{algorithm}[!tb]
\begin{algorithmic}[1]
    \caption{The \(s\)-step GMRES algorithm \label{alg:sstep-GMRES}}
    \REQUIRE
    A matrix \(A \in \mathbb R^{n\times n}\), a right-hand side \(b \in \mathbb R^{n}\), an initial guess \(x^{(0)} \in \mathbb R^{n}\), a block size \(s\), a left-preconditioner \(M_L \in \mathbb R^{n\times n}\), and a right-preconditioner \(M_R \in \mathbb R^{n\times n}\).
    \ENSURE
    A computed solution \(x \in \mathbb R^{n}\) approximating the solution of \(A x = b\).

    \STATE \(r \gets M_L^{-1} (b - A x^{(0)})\) and \(\beta \gets \beta = \norm{r}\).
    \STATE \(V_{1} = r/\beta\). \label{line:sstep-GMRES:V1}
    \FOR{\(i = 1:n/s\)}
        \STATE Perform the \(i\)-th step of the \(s\)-step Arnoldi process (e.g., Algorithm~\ref{alg:sstep-Arnoldi-classical}) to obtain the basis \(B_{1:is}\), the preconditioned basis \(Z_{1:is}\), the orthonormal matrix \(V_{1:is+1}\), and the upper triangular matrix \(R_{1:is+1}\) satisfying \(\bmat{r & W_{1:is}} = V_{1:is+1} R_{1:is+1}\) with \(W_{1:is} \gets M_L^{-1} A Z_{1:is}\).
        \STATE \(H_{1:(i-1)s+1, 1:(i-1)s} \gets R_{1:(i-1)s+1, 2:(i-1)s+1}\).
        \STATE Compute the QR factorization \(H_{1:is+1, (i-1)s+1:is} = G_{1:is+1} T_{1:is+1, (i-1)s+1:is}\) by Givens rotations with an orthogonal matrix \(G_{1:is+1}\), based on \(H_{1:(i-1)s+1, 1:(i-1)s} = G_{1:(i-1)s+1} T_{1:(i-1)s+1, 1:(i-1)s}\).
        \IF{the stopping criterion is satisfied}
             \STATE Solve the triangular system \(T_{1:is} y^{(i)} = \beta G_{1, 1:is}\trans\) to obtain \(y^{(i)} \in \mathbb R^{is}\).
             \RETURN \(x = x^{(i)} \gets x^{(0)} + Z_{1:is} y^{(i)}\).
        \ENDIF
    \ENDFOR
\end{algorithmic}
\end{algorithm}
\section{An abstract framework for backward stability of \(s\)-step GMRES}
\label{sec:stability}

In the recent work~\cite{BHMV2024}, the authors propose a framework to study the backward stability of what they call the ``modular GMRES'' algorithm, which is as follows:
\begin{enumerate}
    \item Compute \(Z_{1:(i+1)s} = M_R^{-1} V_{1:(i+1)s}\) and \(W_{1:(i+1)s} = M_L^{-1} A Z_{1:(i+1)s}\). \label{eq:modular:Z}
    \item Solve \(y^{(i+1)} = \argmin_y \norm{M_L^{-1} b - W_{1:(i+1)s} y}\). \label{eq:modular:ls}
    \item Compute the solution approximation \(x^{(i+1)} = Z_{1:(i+1)s} y^{(i+1)}\). \label{eq:modular:x}
\end{enumerate}
This modular GMRES framework is capable of capturing a wide range of GMRES variants, including the standard GMRES algorithm with various orthogonalization methods.
By analyzing the rounding errors in each step, \cite{BHMV2024} demonstrates that this modular GMRES algorithm is backward stable under mild assumptions.
Among these three steps, Steps~\ref{eq:modular:Z} and~\ref{eq:modular:x} are typically straightforward to analyze using standard rounding error analysis.
To assess the backward stability of GMRES with various orthogonalization methods, it is necessary to evaluate the rounding errors in Step~\ref{eq:modular:ls} for the chosen orthogonalization method within this framework.
Note that there are implicit algorithmic choices in Step~\ref{eq:modular:ls} based on  various orthogonalization methods.
Therefore, Step~\ref{eq:modular:ls} needs to be thoroughly re-evaluated when modifying the orthogonalization method used in GMRES, which is inconvenient (and unnecessary), and further causes ambiguity regarding the requirements of the orthogonalization method.
It is possible to expand the analysis from~\cite{BHMV2024} to the \(s\)-step GMRES method by substituting the orthonormal basis \(V_{1:(i+1)s}\) with a general basis \(B_{1:(i+1)s}\).
However, due to the above-mentioned deficiencies, we will not fully follow the analysis presented in~\cite{BHMV2024}.

In this section, we develop an improved framework for \(Z_{1:(i+1)s} = M_R^{-1} B_{1:(i+1)s}\) focusing on the backward stability of the \(s\)-step GMRES method.
This includes Algorithm~\ref{alg:sstep-GMRES} with the classical $s$-step Arnoldi process (Algorithm~\ref{alg:sstep-Arnoldi-classical}), and the standard GMRES algorithm, which is essentially a specific case of $s$-step GMRES with \(s = 1\) and \(p_0(A) = I\), hence \(B_{1:(i+1)s} = V_{1:(i+1)s}\). Furthermore, we delve into the block orthogonalization, Givens QR decomposition, and triangular system solving steps, which are actually performed in the algorithm, as opposed to analyzing \(y^{(i+1)} = \argmin_y \norm{M_L^{-1} b - W_{1:(i+1)s} y}\).
From this analysis, we will present the assumptions only relevant to the block orthogonalization process itself.
Consequently, we only need to consider the rounding errors in the orthogonalization process when modifying the orthogonalization method utilized in the GMRES algorithm. We note that the approach here thus represents an improved version of the modular GMRES framework in \cite{BHMV2024} for the standard (non $s$-step) GMRES case. 

According to Algorithm~\ref{alg:sstep-GMRES}, the first \((i+1)\) iterations can be summarized as the following four steps.
Without loss of generality, we assume \(x^{(0)} = 0\) for simplicity in the following analysis.
\begin{enumerate}
    \item Compute \(Z_{1:(i+1)s} = M_R^{-1} B_{1:(i+1)s}\) and \(W_{1:(i+1)s} = M_L^{-1} A Z_{1:(i+1)s}\). \label{line:AV}
    \item Compute the QR factorization
    \[
        \bmat{r & W_{1:(i+1)s}} = V_{1:(i+1)s+1} R_{1:(i+1)s+1},
    \]
    where \(V_{1:(i+1)s+1}\) is orthonormal, and \(R_{1:(i+1)s+1}\) is upper triangular. \label{line:orth}
    \item Solve the least squares problem \(y^{(i+1)} = \argmin_y \normbigg{R_{1:(i+1)s+1} \bmat{1 \\ -y}}\), which is through solving the triangular system \(T_{1:(i+1)s} y^{(i+1)} = \beta G_{1, 1:(i+1)s}\trans\) to obtain \(y^{(i+1)}\), where \(T_{1:(i+1)s+1, 1:(i+1)s}\) is the \(R\)-factor of the QR factorization by Givens rotations, i.e.,
    \[
        H_{1:(i+1)s+1, 1:(i+1)s} = G_{1:(i+1)s+1} T_{1:(i+1)s+1, 1:(i+1)s}
    \]
    with \(H_{1:(i+1)s+1, 1:(i+1)s} = R_{1:(i+1)s+1, 2:(i+1)s+1}\) and an orthogonal matrix \(G_{1:(i+1)s+1}\). \label{line:ls}
    \item Compute the solution approximation \(x^{(i+1)} = Z_{1:(i+1)s} y^{(i+1)}\). \label{line:updatex}
\end{enumerate}

Taking rounding errors into account, we further assume that each line satisfies the following, to be described below.
In the remainder of this subsection, \(\varepsilon_*\), \(\delta_*\), and \(\omega_*\) 
denote constants in the rounding error analysis of different operations, which helps to identify the source of errors in the final results.
We use \(\delta\) to indicate constants associated with the left and right preconditioners, and use \(\omega\) to denote constants relevant to the block orthogonalization process.
Constants that are independent of these two parts are denoted by \(\varepsilon\), which can be considered as \(\varepsilon_* \in (0,1)\), as ensured by standard rounding error analysis with the assumption \(n\macheps \ll 1\); see, e.g.,~\cite{H2002}.

\paragraph{Step~\ref{line:AV}: Generating the basis}
This step satisfies
\begin{align}
    \widehat{Z}_{1:(i+1)s} & = M_R^{-1} \widehat B_{1:(i+1)s} + \Delta Z_{1:(i+1)s}, \quad \norm{\Delta Z_{j}} \leq \epsB \norm{M_R^{-1}} \norm{\widehat{B}_{j}}, \label{eq:sstep-GMRES-barB}\\
    \widehat W_{1:(i+1)s} & = M_L^{-1} (A \widehat{Z}_{1:(i+1)s} + \Delta C_{1:(i+1)s}) + \Delta D_{1:(i+1)s}. \label{eq:sstep-GMRES-barW}
\end{align}
Here \(\Delta C_{1:(i+1)s}\) and \(\Delta D_{1:(i+1)s}\) come from, respectively, computing \(A \widehat{Z}_{1:(i+1)s}\) and applying \(M_L^{-1}\) to \(A \widehat{Z}_{1:(i+1)s} + \Delta C_{1:(i+1)s}\), which satisfy
\begin{align*}
    \norm{\Delta C_{j}} & \leq \epsAv \normF{A} \norm{\widehat{Z}_{j}}, \quad
    \norm{\Delta D_{j}} \leq \epsMLAZ \norm{M_L^{-1}} \normF{A} \norm{\widehat{Z}_{j}}
\end{align*}
for any \(j \in \{1, 2, \dotsc, (i+1)s\}\).
Let \(\Delta W_{1:(i+1)s} = M_L^{-1} \Delta C_{1:(i+1)s} + \Delta D_{1:(i+1)s}\).
Then we obtain
\begin{equation} \label{eq:sstep-GMRES-wi}
    \widehat W_{1:(i+1)s} = M_L^{-1} A \widehat{Z}_{1:(i+1)s} + \Delta W_{1:(i+1)s},
\end{equation}
with \(\norm{\Delta W_{j}} \leq (\epsAv + \epsMLAZ) \norm{M_L^{-1}} \normF{A} \norm{\widehat{Z}_{j}}\), which means
\begin{equation} \label{eq:sstep-GMRES-Dwi}
    \normF{\Delta W_{1:(i+1)s} D_{1:(i+1)s}^{-1}} \leq (\epsAv + \epsMLAZ) \norm{M_L^{-1}} \normF{A} \normF{\tilde Z_{1:(i+1)s}}, 
\end{equation}
where \(\widehat Z_{1:(i+1)s} = \tilde Z_{1:(i+1)s} D_{1:(i+1)s}\) with any invertible diagonal matrix \(D_{1:(i+1)s}\).

\paragraph{Step~\ref{line:orth}: Backward stability of the block orthogonalization}
Let \(r = M_L^{-1} (b - A x^{(0)}) = M_L^{-1} b\) with \(x^{(0)} = 0\) and \(\beta = \norm{\widehat{r}}\) satisfy
\begin{align}
    & \widehat{r} = M_L^{-1} b + \Delta r, \quad \norm{\Delta r} \leq \epsMLb \norm{M_L^{-1}} \norm{b},  \label{eq:sstep-GMRES:barr} \\
    & \widehat \beta = (1 + \epsr) \beta. \label{eq:sstep-GMRES:beta}
\end{align}
Notice that the first column of the \(Q\)-factor associated with \(\bmat{\phi \widehat{r}& \widehat{W}_{1:(i+1)s}}\), specifically \(V_{1:(i+1)s+1}\), is determined by Line~\ref{line:sstep-GMRES:V1} in Algorithm~\ref{alg:sstep-Arnoldi-classical}, and this computation is column-wise backward stable.
Thus, in GMRES, the orthogonalization method employed in Step~\ref{line:orth} can be viewed as being applied to the matrix \(\bmat{\phi \widehat{r}& \widehat{W}_{1:(i+1)s}}\) for any positive scalar \(\phi\).
We assume the backward stability of the block orthogonalization method used in Step~\ref{line:orth}; 
according to~\cite{CLRT2022}, most commonly-used block orthogonalization schemes satisfy this property,
i.e., for any \(\phi > 0\),
\begin{equation} \label{eq:sstep-GMRES-VR}
    \bmat{\phi \widehat{r}& \widehat{W}_{1:(i+1)s}} + \Delta E_{1:(i+1)s+1}(\phi) = \bar{V}_{1:(i+1)s+1} \widehat{R}_{1:(i+1)s+1}(\phi),
\end{equation}
where \(\widehat{R}_{1:(i+1)s+1}(\phi) = \bmat{\widehat{R}_{1:(i+1)s+1, 1}(\phi) & \widehat{R}_{1:(i+1)s+1, 2:(i+1)s+1}}\)
and the backward error \(\Delta E_{1:(i+1)s+1}(\phi) = \bmat{\Delta E_1(\phi) & \Delta E_{2:(i+1)s+1}}\) satisfies
\begin{align}
    & \norm{\Delta E_1(\phi)} \leq \epsqr \norm{\phi \widehat{r}} = \epsqr \phi \beta \quad \text{and} 
    \quad \norm{\Delta E_j} \leq \epsqr \norm{\widehat W_{j-1}} \label{eq:sstep-GMRES-VR-err}
\end{align}
for any \(j \in \{2, \dotsc, (i+1)s+1\}\).
Furthermore, \eqref{eq:sstep-GMRES-VR-err} implies
\begin{equation} \label{eq:sstep-GMRES-DEyi}
\begin{split}
    & \begin{split}
        & \normFbigg{\Delta E_{1:is+j}(\phi) \bmat{1 & 0 \\ 0 & D_{1:is+j-1}^{-1}}} \\
        & \quad\leq \epsqr \phi \beta + \epsqr \normF{\widehat{W}_{1:is+j} D_{1:is+j-1}^{-1}} \\
        & \quad \leq \epsqr \phi \beta + \epsqr \bigl(1 + \epsAv + \epsMLAZ\bigr) \norm{M_L^{-1}} \normF{A} \normF{\tilde Z_{is+j-1}}, \\
    \end{split} \\
    & \begin{split}
        &\normbigg{\Delta E_{1:(i+1)s+1}(\phi) \bmat{1\\ -\widehat{y}^{(i+1)}}} \\
        & \quad \leq \normbigg{\Delta E_{1:(i+1)s+1}(\phi) \bmat{1 & 0 \\ 0 & D_{1:(i+1)s}^{-1}} \bmat{1\\ -D_{1:(i+1)s} \widehat{y}^{(i+1)}}} \\
        & \quad \leq \epsqr \phi \beta + \epsqr \bigl(1 + \epsAv + \epsMLAZ\bigr) \norm{M_L^{-1}} \normF{A} \normF{\tilde Z_{1:(i+1)s}} \norm{D_{1:(i+1)s} \widehat{y}^{(i+1)}},
    \end{split}
\end{split}
\end{equation}
for any \(j \in \{1, 2, \dotsc, s\}\).
Here \(\bar V_{*}\) can be the computed result \(\widehat V_{*}\) if using the reorthogonalized block classical Gram--Schmidt algorithm, or an exact orthonormal matrix \(\tilde V_{*}\) if using the block modified Gram--Schmidt algorithm.

In addition, we define \(\omega_{is+j}\) to describe the loss of orthogonality of \(\bar V_{1:is+j}\), i.e., for any \(j \in \{1, 2, \dotsc, s\}\),
\begin{equation} \label{eq:sstep-GMRES-Viorth}
\begin{split}
    \normF{\bar V_{1:is+j}\trans \bar V_{1:is+j} - I} &\leq \omega_{is+j}, \\
    \norm{\bar V_{1:is+j}} &= \sqrt{\norm{\bar V_{1:is+j}\trans \bar V_{1:is+j}}} \leq \sqrt{\norm{I} + \norm{\bar V_{1:is+j}\trans \bar V_{1:is+j} - I}} \\
    &\leq \sqrt{1 + \omega_{is+j}}.
\end{split}
\end{equation}

\paragraph{Step~\ref{line:orth}: Loss of orthogonality in the block orthogonalization}
In this part, we discuss the iteration where \(\widehat V_*\) is no longer well-conditioned, which is equivalent to the condition that the columns of \(\widehat V_*\) have lost orthogonality.
This iteration is called the ``key dimension'' in~\cite{BHMV2024}, since the loss of orthogonality always happens when the condition number of \(\bmat{\phi r & \widehat{W}_{1:p} D_{1:p}^{-1}}\) becomes large enough. 
We employ a similar definition of the key dimension introduced by~\cite{BHMV2024} here. We define \(p = \ki s + \kj\) to be the iteration in which we reach the key dimension, i.e., when we have 
\begin{align}
    & \sigmin\Big(\bmat{\phi \widehat{r}& \widehat{W}_{1:p} D_{1:p}^{-1}}\Big) \leq \epskapparW \normFBig{\bmat{\phi \widehat{r}& \widehat{W}_{1:p} D_{1:p}^{-1}}}. \label{eq:lem:LS:assump-rW}
\end{align}
Together with
\begin{equation}
    \sigmin(\widehat{W}_{1:p} D_{1:p}^{-1}) \geq \epskappaW \normF{\widehat{W}_{1:p} D_{1:p}^{-1}},  \label{eq:lem:LS:assump-W}
\end{equation}
\eqref{eq:lem:LS:assump-rW} indicates that \(\widehat{r}\) lies in the range of \(\widehat W_{1:p}\).

\paragraph{Step~\ref{line:ls}: Solving the least squares problem}

We only need to consider the case of the key dimension.
The procedure of solving the least squares problem in Step~\ref{line:ls} at the key dimension satisfies
\begin{align}
    \widehat y^{(k)} & = \argmin_y \normbigg{(\widehat{R}_{1:p+1} + \Delta R_{1:p+1}) \bmat{1 \\ -y}}, \quad
    \norm{\Delta R_{1:p+1, j}} \leq \epsls \norm{\widehat R_{1:p+1, j}}, \label{eq:sstep-GMRES-ls}
\end{align}
for any \(j = 1, \dotsc, p+1\).
The error \(\Delta R_{1:p+1}\) further satisfies
\begin{equation} \label{eq:sstep-GMRES-ls-er}
    \normFBig{\bmat{\phi \Delta R_{1:p+1, 1} & \Delta R_{1:p+1, 2:p+1} D_{1:p}^{-1}}} \leq \epsls \normFBig{\bmat{\phi \widehat R_{1:p+1, 1} & \widehat R_{1:p+1, 2:p+1} D_{1:p}^{-1}}}.
\end{equation}

Regarding the least squares problem, \cite{BHMV2024} measures \(\min_y {\norm{\widehat r - \widehat W_{1:p} y}}\), which involves Steps~\ref{line:orth} and~\ref{line:ls} of the algorithm.
Thus, the substantial impact of the block orthogonalization process, i.e., Step~\ref{line:orth}, will be unclear when assessing \(\min_y {\norm{\widehat r - \widehat W_{1:p} y}}\) directly.
In this context, we avoid directly accessing \(\min_y {\norm{\widehat r - \widehat W_{1:p} y}}\) and aim to distinguish the errors introduced by Steps~\ref{line:orth} and~\ref{line:ls} separately.
Therefore, our attention is on~\eqref{eq:sstep-GMRES-ls}, which is performed in Step~\ref{line:ls} and remains unaffected by the block orthogonalization.

\paragraph{Step~\ref{line:updatex}: Updating \(x\)}
The step of updating \(x\) satisfies
\begin{equation}
    \widehat{x}^{(\ki)} = \widehat{Z}_{1:p} \widehat{y}^{(\ki)} + \Delta x^{(\ki)}, \quad \norm{\Delta x^{(\ki)}} \leq \epsMV \normF{\tilde Z_{1:p}} \norm{D_{1:p} \widehat{y}^{(\ki)}}. \label{eq:sstep-GMRES-xi}
\end{equation}

We summarize the quantities \(\epsilon_{*}\), \(\delta_*\), and \(\omega_*\), coming from different operations in Table~\ref{table:symbols}.
\begin{table}[!tb]
\centering
\caption{Notation of rounding error analysis for Steps~\ref{line:AV}, \ref{line:orth}, \ref{line:ls}, and~\ref{line:updatex}.}
\begin{tabular}{ccc}
\hline
Notation          & Sources                                                           & Details                                        \\
\hline
\(\epsB\)         & Applying \(M_R^{-1}\) to \(B_{1:(i+1)s}\)                                       & \eqref{eq:sstep-GMRES-barB}   \\
\(\epsAv\)        & Computing \(AZ_{1:(i+1)s}\)                                                     & \eqref{eq:sstep-GMRES-barW}   \\
\(\epsMLAZ\)      & Applying \(M_L^{-1}\) to \(AZ_{1:(i+1)s}\)                                      & \eqref{eq:sstep-GMRES-barW}   \\
\(\epsMLb\)       & Applying \(M_L^{-1}\) to \(b\)                                                  & \eqref{eq:sstep-GMRES:barr}   \\
\(\epsqr\)        & Orthogonalizing \(\bmat{\phi r& W_{1:(i+1)s}}\)                                 & \eqref{eq:sstep-GMRES-VR}     \\
\(\omega_{is+j}\) & Loss of orthogonality of \(V_{1:is+j}\)                                         & \eqref{eq:sstep-GMRES-Viorth} \\
\(\epskapparW\)   & \makecell{The upper bound on \(1/\kappa\bigl(\bmat{\phi \widehat{r}& \widehat{W}_{1:p} D_{1:p}^{-1}}\bigr)\) to show \\ if \(r\) lies in the range of \(W_{1:p}\)}                                          & \eqref{eq:lem:LS:assump-rW}   \\
\(\epskappaW\)    & \makecell{The lower bound on \(1/\kappa(\widehat{W}_{1:p} D_{1:p}^{-1})\) to show \\ if the columns of \(W_{1:p}\) are linearly independent}                             & \eqref{eq:lem:LS:assump-W}    \\
\(\epsls\)        & Solving the least squares problem \(\min_y \normbigg{R_{1:p+1}\bmat{1 \\ -y}}\) & \eqref{eq:sstep-GMRES-ls}     \\
\(\epsMV\)        & Computing \(Z_{1:p}y^{(\ki)}\)                                                  & \eqref{eq:sstep-GMRES-xi}   \\
\hline
\end{tabular}
\label{table:symbols}
\end{table}

\begin{remark} \label{remark:eps-1}
    From standard rounding error analysis and~\cite[Equation~(3.12)]{H2002}, we have \(\epsr\), \(\epsAv\), \(\epsMV = \bigO(\macheps)\).
    In~\cite[Equation~(A.7)]{BHMV2024}, the authors have given \(\epsls = \bigO(\macheps)\), 
    which comes from applying Givens QR factorization to the upper Hessenberg matrix \(\hat R_{1:p+1, 2:p+1}\) and solving the triangular system as depicted in Step~\ref{line:ls}.

    The terms \(\epsMLAZ\), \(\epsMLb\), and \(\epsB\) are highly dependent on the preconditioners \(M_L\) and \(M_R\).
    The remaining terms \(\epsqr\), \(\omega_{is+j}\), \(\epskapparW\), and \(\epskappaW\) are determined by the block orthogonalization algorithm employed in Step~\ref{line:orth}.
\end{remark}

In the implementation of the $s$-step GMRES algorithm, the stopping criterion typically involves checking whether \(\min_y \normbigg{(\widehat{R}_{1:is+j+1} + \Delta R_{1:is+j+1}) \bmat{1 \\ -y}}\) is sufficiently small.
Thus, we first show that the residual \(\normbigg{\widehat{R}_{1:p+1} \bmat{1 \\ -\widehat y^{(k)}}}\) is small enough in Lemma~\ref{lem:LS}.
Then we prove that the backward error \(\norm{b - A \widehat x^{(k)}}\) can be bounded mainly by the residual of the least squares problem in Theorem~\ref{thm:b-Axi-H}.

According to~\cite[Theorem 2.4]{PRS2006}, the residual of the least squares problem in Step~\ref{line:ls}, i.e., \(\normbigg{\widehat{R}_{1:p+1} \bmat{1 \\ -\widehat y^{(k)}}}\), can be bounded by \(\sigmin\bigg(\bmat{\phi \widehat R_{1:p+1, 1} & \widehat R_{1:p+1, 2:p+1} D_{1:p}^{-1}}\bigg)\) and \(\normFBig{\bmat{\phi \widehat R_{1:p+1, 1} & \widehat R_{1:p+1, 2:p+1} D_{1:p}^{-1}}}\).
Therefore, establishing a bound for the residual \(\normbigg{\widehat{R}_{1:p+1} \bmat{1 \\ -\widehat y^{(k)}}}\) requires determining the relationship between the result and the assumptions~\eqref{eq:lem:LS:assump-rW} and \eqref{eq:lem:LS:assump-W}, which amounts to bounding 
\[
    \sigmin\bigg(\bmat{\phi \widehat R_{1:p+1, 1} & \widehat R_{1:p+1, 2:p+1} D_{1:p}^{-1}}\bigg)
\]
by \(\normFBig{\bmat{\phi r, \widehat{W}_{1:p} D_{1:p}^{-1}}}\) and bounding 
\[
    \normFBig{\bmat{\phi \widehat R_{1:p+1, 1} & \widehat R_{1:p+1, 2:p+1} D_{1:p}^{-1}}}
\]
by \(\normF{\widehat{W}_{1:p} D_{1:p}^{-1}}\).
We summarize the result in the following lemma.

\begin{lemma} \label{lem:LS-norm}
    Assume that \(\bmat{\phi r, \widehat{W}_{1:p} D_{1:p}^{-1}}\) and \(\widehat{W}_{1:p} D_{1:p}^{-1}\) satisfy~\eqref{eq:lem:LS:assump-rW} and~\eqref{eq:lem:LS:assump-W}, respectively.
    If also assuming \(\normF{\bar V_{1:p}\trans \bar V_{1:p} - I} \leq \omega_{p} \leq 1/2\), then
    \begin{align}
        \sigmin\Big(\bmat{\phi \widehat R_{1:p+1, 1}, \widehat R_{1:p+1, 2:p+1} D_{1:p}^{-1}}\Big) \leq{}& \alpha_{11} \normFBig{\bmat{\phi r, \widehat{W}_{1:p} D_{1:p}^{-1}}} + \alpha_{12} \normF{\widehat W_{1:p} D_{1:p}^{-1}}, \label{eq:lem-LS-norm:sigminbigR-sum} \\
        \normFBig{\bmat{\phi \widehat R_{1:p+1, 1}, \widehat R_{1:p+1, 2:p+1} D_{1:p}^{-1}}} \leq{}& \alpha_{21} \normFBig{\bmat{\phi r, \widehat{W}_{1:p} D_{1:p}^{-1}}} + \alpha_{22} \normF{\widehat W_{1:p} D_{1:p}^{-1}}, \label{eq:lem-LS-norm:normbigR-sum}\\
        \sigmin(\widehat R_{1:p+1, 2:p+1} D_{1:p}^{-1}) \geq{}& \alpha_3 \normF{\widehat W_{1:p} D_{1:p}^{-1}}, \label{eq:lem-LS-norm:sigminR-sum} \\
        \normF{\widehat R_{1:p+1, 2:p+1} D_{1:p}^{-1}} \leq{}& \alpha_4 \normF{\widehat W_{1:p} D_{1:p}^{-1}} \label{eq:lem-LS-norm:normR-sum},
    \end{align}
    where
    \begin{align}
        & \alpha_{11} = 3 (\epskapparW + \epsqr), \quad \alpha_{12} = 4 \cdot \epsH, \\
        & \alpha_{21} = 3 (1 + \epsqr), \quad \alpha_{22} = 4 \cdot \epsH, \\
        & \alpha_3 = \frac{2}{3} (\epskappaW - \epsqr - \epsH) - \epsH, \\
        & \alpha_4 = 3 (1 + \epsqr + \epsH) + \epsH,
    \end{align}
    with \(\epsH\) defined by
    \begin{equation} \label{eq:lem-LS:assump-epsH}
        \max\bigl(\abs{\widehat R_{p+1, p+1} D_{p, p}^{-1}}, \normF{\bar V_{p+1} \widehat R_{p+1, 2:p+1} D_{1:p}^{-1}}\bigr) \leq \epsH \normF{\widehat W_{1:p} D_{1:p}^{-1}},
    \end{equation}
    if \(\normF{\bar V_{1:p}\trans \bar V_{1:p} - I} > 0\); otherwise, \(\epsH = 0\).
\end{lemma}

Notice that \(\hat R_{1:p+1}\) is obtained as the \(R\)-factor by employing the block QR factorization on \(\bmat{\phi \widehat{r}& \widehat{W}_{1:(i+1)s}}\) as indicated in~\eqref{eq:sstep-GMRES-VR}.
Therefore, before proving Lemma~\ref{lem:LS-norm}, we present the following lemma to estimate the largest and smallest singular values of the computed \(R\) factor.

\begin{lemma} \label{lem:sigmaR}
    Given \(X \in \mathbb{R}^{n \times m}\),
    assume that \(X + \Delta X = Q U\), where \(Q \in \mathbb{R}^{n \times m}\) satisfies \(\normF{Q\trans Q - I} < 1\).
    Then
    \begin{align}
        & \normF{U} \leq \frac{\norm{Q} (\normF{X} + \normF{\Delta X})}{1 - \normF{Q\trans Q - I}} \quad \text{and} \quad
        \sigmin(U) \geq \frac{\sigmin(X) - \normF{\Delta X}}{\norm{Q}}.
    \end{align}
\end{lemma}

\begin{proof}
    By \(X + \Delta X = Q U\), we have \(U = Q\trans X + Q\trans \Delta X - (Q\trans Q - I) U\) and then
    \begin{equation}
        \normF{U} \leq \norm{Q} (\normF{X} + \normF{\Delta X}) + \norm{Q\trans Q - I} \normF{U},
    \end{equation}
    which gives the bound of \(\normF{U}\).
    The bound of \(\sigmin(U)\) has been proved in~\cite[Lemma 2]{CLMO2024}.
\end{proof}

We now prove Lemma~\ref{lem:LS-norm}.

\begin{proof}[Proof of Lemma~\ref{lem:LS-norm}]
    Notice that \(\bmat{\phi \widehat R_{1:p+1, 1} & \widehat R_{1:p+1, 2:p+1}}\) is the computed \(R\)-factor of \(\bmat{\phi r, \widehat{W}_{1:p}}\) from~\eqref{eq:sstep-GMRES-VR}.
    We then multiply two sides of ~\eqref{eq:sstep-GMRES-VR} by \(\bmat{1 & 0\\ 0& D^{-1}_{1:p}}\) to obtain
    \begin{equation} \label{eq:sstep-GMRES-VRD}
        \bmat{\phi \widehat{r}& \widehat{W}_{1:p}D_{1:p}^{-1}} + \Delta \tilde E_{1:p+1}(\phi) = \bar{V}_{1:p+1} \bmat{\phi \widehat R_{1:p+1, 1} & \widehat R_{1:p+1, 2:p+1} D_{1:p}^{-1}},
    \end{equation}
    where \(\Delta \tilde E_{1:p+1}(\phi):= \Delta E_{1:p+1}(\phi)\bmat{1 & 0\\ 0& D^{-1}_{1:p}}\) satisfies~\eqref{eq:sstep-GMRES-DEyi}.
    
    From~\eqref{eq:sstep-GMRES-VRD}, if \(\bar V_{1:p+1}\) is exactly orthonormal, we have
    \begin{align*}
        \sigmin\Big(\bmat{\phi \widehat R_{1:p+1, 1} & \widehat R_{1:p+1, 2:p+1} D_{1:p}^{-1}}\Big) &\leq \sigmin\Big(\bmat{\phi r, \widehat{W}_{1:p} D_{1:p}^{-1}}\Big) \\
        &\quad + \epsqr \normFBig{\bmat{\phi r, \widehat{W}_{1:p} D_{1:p}^{-1}}}, \\
        \normFBig{\bmat{\phi \widehat R_{1:p+1, 1} & \widehat R_{1:p+1, 2:p+1} D_{1:p}^{-1}}} &\leq (1 + \epsqr) \normFBig{\bmat{\phi r, \widehat{W}_{1:p} D_{1:p}^{-1}}}, \\
        \sigmin(\widehat R_{1:p+1, 2:p+1} D_{1:p}^{-1}) &\geq \sigmin(\widehat{W}_{1:p} D_{1:p}^{-1}) - \epsqr \normF{\widehat{W}_{1:p} D_{1:p}^{-1}}, \\
        \normF{\widehat R_{1:p+1, 2:p+1} D_{1:p}^{-1}} &\leq (1 + \epsqr) \normF{\widehat{W}_{1:p} D_{1:p}^{-1}},
    \end{align*}
    which proves the result under the assumptions~\eqref{eq:lem:LS:assump-rW} and~\eqref{eq:lem:LS:assump-W}.

    Then we consider the case when \(\bar V_{1:p+1}\) is close orthonormal instead of being exactly orthonormal, i.e, \(0 < \normF{\bar V_{1:p}\trans \bar V_{1:p} - I} \leq \omega_{p} \leq 1/2\).
    In this case, we need to avoid to directly considering \(\bar V_{1:p+1}\) since \(\bar V_{1:p+1}\) may be not well-conditioned.
    Observe that
    \[
        \bar V_{1:p+1} \widehat R_{1:p+1, 2:p+1} = \bar V_{1:p} \widehat R_{1:p, 2:p+1} + \bar V_{p+1} \widehat R_{p+1, 2:p+1}
    \]
    and further
    \begin{equation} \label{eq:lemma-LS-norm:VR-split}
        \begin{split}
            \bar V_{1:p+1} & \bmat{\phi \widehat R_{1:p+1, 1} & \widehat R_{1:p+1, 2:p+1} D_{1:p}^{-1}} \\
            &= \bar V_{1:p} \bmat{\phi \widehat R_{1:p, 1} & \widehat R_{1:p, 2:p+1} D_{1:p}^{-1}} + \bar V_{p+1} \widehat R_{p+1, 1:p+1} D_{1:p}^{-1}.
        \end{split}
    \end{equation}
    By substituting~\eqref{eq:lemma-LS-norm:VR-split} into the right-hand side of~\eqref{eq:sstep-GMRES-VRD}, we derive
    \begin{equation} \label{eq:lemma-LS-norm:QR-rW}
        \begin{split}
            [\phi r, \widehat W_{1:p} D_{1:p}^{-1}] + \Delta F = \bar V_{1:p} \bmat{\phi \widehat R_{1:p, 1} & \widehat R_{1:p, 2:p+1} D_{1:p}^{-1}}
        \end{split}
    \end{equation}
    with \(\Delta F = \Delta\tilde E_{1:p+1}(\phi) - \bar V_{p+1} \widehat R_{p+1, 1:p+1} D_{1:p}^{-1}\) satisfying, from~\eqref{eq:sstep-GMRES-DEyi} and the definition~\eqref{eq:lem-LS:assump-epsH} of \(\epsH\),
    \[
        \normF{\Delta F} \leq \epsqr \normFbig{[\phi r, \widehat W_{1:p} D_{1:p}^{-1}]} + \epsH \normF{\widehat W_{1:p} D_{1:p}^{-1}}.
    \]
    In~\eqref{eq:lemma-LS-norm:QR-rW}, notice that \(\bmat{\phi \widehat R_{1:p, 1} & \widehat R_{1:p, 2:p+1} D_{1:p}^{-1}}\) is the \(R\)-factor of \([\phi r, \widehat W_{1:p} D_{1:p}^{-1}] + \Delta F\).
    Using~\eqref{eq:lem:LS:assump-rW} and the perturbation results of QR factorization, i.e., Lemma~\ref{lem:sigmaR}, we have the bound on the largest singular value,
    \begin{equation}
        \begin{split}
            &\normbig{\bmat{\phi \widehat R_{1:p, 1} & \widehat R_{1:p, 2:p+1} D_{1:p}^{-1}]}} \\
            & \qquad \leq \frac{1 + \omega_p}{1 - \omega_p} \bigl( (1 + \epsqr) \normFbig{[\phi r, \widehat W_{1:p} D_{1:p}^{-1}]} + \epsH \normFbig{\widehat W_{1:p} D_{1:p}^{-1}} \bigr) \\
            & \qquad \leq 3 \bigl( (1 + \epsqr) \normFbig{[\phi r, \widehat W_{1:p} D_{1:p}^{-1}]} + \epsH \normFbig{\widehat W_{1:p} D_{1:p}^{-1}} \bigr),
        \end{split}
    \end{equation}
    and the smallest one,
    \begin{equation}
        \begin{split}
            \sigmin\big([\phi \widehat R_{1:p, 1} & \quad \widehat R_{1:p, 2:p+1} D_{1:p}^{-1}]\big) \\
            \leq{}& 3 \bigl( (\epskapparW + \epsqr) \normFbig{[\phi r, \widehat W_{1:p} D_{1:p}^{-1}]} + \epsH \normFbig{\widehat W_{1:p} D_{1:p}^{-1}} \bigr).
        \end{split}
    \end{equation}
    By substituting these two bounds and~\eqref{eq:lem-LS:assump-epsH} into the right-hand side of the two inequalities
    \[
       \normbig{[\phi \widehat R_{1:p+1, 1} \quad \widehat R_{1:p+1, 2:p+1} D_{1:p}^{-1}]}
        \leq \normbig{[\phi \widehat R_{1:p, 1} \quad \widehat R_{1:p, 2:p+1} D_{1:p}^{-1}]} + \abs{\widehat R_{p+1, p+1} D_{p, p}^{-1}}
    \]
    and
    \begin{equation*}
        \begin{split}
            & \sigmin\big([\phi \widehat R_{1:p+1, 1} \quad \widehat R_{1:p+1, 2:p+1} D_{1:p}^{-1}]\big) \\
            & \qquad \leq \sigmin\big([\phi \widehat R_{1:p, 1} \quad \widehat R_{1:p, 2:p+1} D_{1:p}^{-1}]\big) + \abs{\widehat R_{p+1, p+1} D_{p, p}^{-1}},
        \end{split}
    \end{equation*}
    we prove~\eqref{eq:lem-LS-norm:sigminbigR-sum} and~\eqref{eq:lem-LS-norm:normbigR-sum}.

    It remains to prove~\eqref{eq:lem-LS-norm:sigminR-sum} and~\eqref{eq:lem-LS-norm:normR-sum}, which is similar to the proof of~\eqref{eq:lem-LS-norm:sigminbigR-sum} and~\eqref{eq:lem-LS-norm:normbigR-sum}.
    From Lemma~\ref{lem:sigmaR} and
    \begin{align*}
        & \widehat R_{1:p+1, 2:p+1} D_{1:p}^{-1} = \widehat R_{1:p, 2:p+1} D_{1:p}^{-1} + \widehat R_{p+1, 2:p+1} D_{1:p}^{-1}, \\
        & \widehat W_{1:p} D_{1:p}^{-1} + \Delta E_{2:p+1}(\phi) D_{1:p}^{-1} - \bar V_{p+1} \widehat R_{p+1, 2:p+1} D_{1:p}^{-1} = \bar V_{1:p} \widehat R_{1:p, 2:p+1} D_{1:p}^{-1},
    \end{align*}
    we have
    \begin{equation*}
        \begin{split}
            \sigmin(\widehat R_{1:p+1, 2:p+1} D_{1:p}^{-1})  \geq{}& \sigmin(\widehat R_{1:p, 2:p+1} D_{1:p}^{-1}) - \abs{\widehat R_{p+1, p+1} D_{p, p}^{-1}} \\
            \geq{}& \frac{1}{1 + \omega_p} \bigl(\sigmin(\widehat W_{1:p} D_{1:p}^{-1}) - \normF{\Delta E_{2:p+1}(\phi) D_{1:p}^{-1}} \\
            & - \normFbig{\bar V_{p+1} \widehat R_{p+1, 2:p+1} D_{1:p}^{-1}}\bigr) - \abs{\widehat R_{p+1, p+1} D_{p, p}^{-1}} \\
            \geq{}& \biggl(\frac{2}{3} (\epskappaW - \epsqr - \epsH) - \epsH\biggr) \normFbig{\widehat W_{1:p} D_{1:p}^{-1}},
        \end{split}
    \end{equation*}
    and similarly,
    \begin{equation*}
        \normFbig{\widehat R_{1:p+1, 2:p+1} D_{1:p}^{-1}} \leq \biggl(\frac{1 + \omega_p}{1 - \omega_p} (1 + \epsqr + \epsH) + \epsH\biggr) \normFbig{\widehat W_{1:p} D_{1:p}^{-1}},
    \end{equation*}
    which proves~\eqref{eq:lem-LS-norm:sigminR-sum} and~\eqref{eq:lem-LS-norm:normR-sum}.
\end{proof}

Lemma~\ref{lem:LS-norm} gives bounds on the smallest and largest singular values of the matrices \(\bmat{\phi \widehat R_{1:p+1, 1}& \widehat R_{1:p+1, 2:p+1} D_{1:p}^{-1}}\) and \(\widehat R_{1:p+1, 2:p+1} D_{1:p}^{-1}\), which will be required for the proof of Lemma~\ref{lem:LS}.
We are now prepared to bound the residual of the least squares problem in Step~\ref{line:ls}, i.e., \(\normbigg{\widehat{R}_{1:p+1} \bmat{1 \\ -\widehat y^{(k)}}}\), using~\cite[Theorem 2.4]{PRS2006}.

\begin{lemma} \label{lem:LS}
    Assume that \(\widehat y^{(k)}\) satisfies~\eqref{eq:sstep-GMRES-ls}.
    If \eqref{eq:lem:LS:assump-rW} and~\eqref{eq:lem:LS:assump-W} hold with
    \begin{equation} \label{eq:lem-LS:assump}
        \begin{split}
            \epskappaW
            \geq{}& \frac{27 \bigl(\epskapparW + \epsqr + \epsls \bigr)}
            {1 - \epsr - 6 \cdot (\epskapparW + \epsqr) - 9 \cdot \epsls}
            + 9\cdot \epskapparW \\
            & + 10 \cdot \epsqr + 9 \cdot \epsls + 16 \cdot \epsH,
    \end{split}
    \end{equation}
    then
    \begin{equation}
        \begin{split}
            \normbigg{\widehat{R}_{1:p+1} \bmat{1 \\ -\widehat y^{(k)}}}
            \leq{}& (9 \cdot \epskapparW + 9 \cdot \epsqr + 12 \cdot \epsls) \beta
            + (9 \cdot \epskapparW + 9 \cdot \epsqr \\
            & + 12 \cdot \epsls + 12 \cdot \epsH) \normF{\widehat{W}_{1:p} D_{1:p}^{-1}} \norm{D_{1:p} \widehat{y}^{(k)}},
        \end{split}
    \end{equation}
    where \(\epsH\) is defined in Lemma~\ref{lem:LS-norm}.
\end{lemma}

\begin{proof}
    We follow the approach from~\cite{BHMV2024,PRS2006} and employ~\cite[Theorem 2.4]{PRS2006} to estimate the residual of the least squares problem.
    Let \(\tilde R_{1:p+1} = \widehat R_{1:p+1} + \Delta R_{1:p+1}\).
    Then using~\cite[Theorem 2.4]{PRS2006}, the residual of the least squares problem can be estimated as, for any \(\phi > 0\),
    \begin{equation} \label{eq:lem-LS:norm-Ry-1}
        \begin{split}
            \normbigg{(\widehat{R}_{1:p+1} & + \Delta R_{1:p+1}) \bmat{1 \\ -\widehat y^{(k)}}}^2 \\
            \leq{}& \min_y \normbigg{(\widehat{R}_{1:p+1} + \Delta R_{1:p+1}) \bmat{1 & 0 \\ 0 & D_{1:p}^{-1}} \bmat{1 \\ -D_{1:p} y}}^2 \\
            \leq{}& \sigmin^2\big(\bmat{\phi \tilde R_{1:p+1, 1} & \tilde R_{1:p+1, 2:p+1} D_{1:p}^{-1}}\big) \biggl(\frac{1}{\phi^2} + \frac{\normbig{D_{1:p} \widehat{y}^{(i)}}^2}{1 - \delta^2(\phi)}\biggr),
        \end{split}
    \end{equation}
    where \(\delta(\phi)\) is defined by
    \begin{equation} \label{eq:lem-LS:defdelta}
    \delta(\phi) = \frac{\sigmin\big(\bmat{\phi \tilde R_{1:p+1, 1} & \tilde R_{1:p+1, 2:p+1} D_{1:p}^{-1}}\big)}{\sigmin(\tilde R_{1:p+1, 2:p+1} D_{1:p}^{-1})}.
    \end{equation}
    It is clear that \(\delta(\phi)<1\) for any \(\phi>0\).
    Analogous to~\cite{BHMV2024,PRS2006}, we choose
    \begin{equation} \label{eq:lem-LS:defphi}
        \frac{1}{\phi^2} = \frac{\normbig{D_{1:p} \widehat{y}^{(i)}}^2}{1 - \delta^2(\phi)}
    \end{equation}
    to simplify~\eqref{eq:lem-LS:norm-Ry-1} as
    \begin{equation} \label{eq:lem-LS:norm-Ry}
        \begin{split}
            \normbigg{(\widehat{R}_{1:p+1} & + \Delta R_{1:p+1}) \bmat{1 \\ -\widehat y^{(k)}}}
            \leq \sqrt{2}\,\sigmin\big(\bmat{\phi \tilde R_{1:p+1, 1} & \tilde R_{1:p+1, 2:p+1} D_{1:p}^{-1}}\big) \phi^{-1}.
        \end{split}
    \end{equation}
    From~\cite[Equations~(3.17)--(3.18)]{BHMV2024}, we can similarly check that there exists a \(0<\phi<\normbig{D_{1:p} \widehat{y}^{(i)}}^{-1}\) satisfying~\eqref{eq:lem-LS:defphi} and \(\delta(\phi)<1\).
    Then our aim is to bound, respectively, \(\sigmin\big(\bmat{\phi \tilde R_{1:p+1, 1} & \tilde R_{1:p+1, 2:p+1} D_{1:p}^{-1}}\big)\) and \(\delta(\phi)\).
    
    First, we will give the bound for \(\sigmin\big(\bmat{\phi \tilde R_{1:p+1, 1} & \tilde R_{1:p+1, 2:p+1} D_{1:p}^{-1}}\big)\).
    Recalling the definition of \(\tilde R_{1:p+1}\), i.e., \(\tilde R_{1:p+1} = \widehat R_{1:p+1} + \Delta R_{1:p+1}\), we bound the smallest singular value of \([\phi \tilde R_{1:p+1, 1} \quad \tilde R_{1:p+1, 2:p+1} D_{1:p}^{-1}]\) from the perturbation theory of singular values as follows:
    \begin{equation}  \label{eq:lem-LS:sigminbigR-1}
        \begin{split}
            &\sigmin\big([\phi \tilde R_{1:p+1, 1} \quad \tilde R_{1:p+1, 2:p+1} D_{1:p}^{-1}]\big) \\
            &\leq \sigmin\Big(\bmat{\phi \widehat R_{1:p+1, 1} & \widehat R_{1:p+1, 2:p+1} D_{1:p}^{-1}}\Big)
            + \normFBig{\bmat{\phi \Delta R_{1:p+1, 1} & \Delta R_{1:p+1, 2:p+1} D_{1:p}^{-1}}}.
        \end{split}
    \end{equation}
    Using~\eqref{eq:lem-LS-norm:sigminbigR-sum} in Lemma~\ref{lem:LS-norm} and~\eqref{eq:sstep-GMRES-ls-er}, we can bound the first and the second terms of the bound in~\eqref{eq:lem-LS:sigminbigR-1}, respectively, i.e.,
    \begin{equation} \label{eq:lem-LS:sigminbigR-2}
        \begin{split}
            & \sigmin\big([\phi \tilde R_{1:p+1, 1} \quad \tilde R_{1:p+1, 2:p+1} D_{1:p}^{-1}]\big) \\
            \leq{}& \alpha_{11}\normFBig{\bmat{\phi r, \widehat{W}_{1:p} D_{1:p}^{-1}}}
            + \alpha_{12}\normFbig{\widehat{W}_{1:p} D_{1:p}^{-1}} + \epsls \normFBig{\bmat{\phi \widehat R_{1:p+1, 1} & \widehat R_{1:p+1, 2:p+1} D_{1:p}^{-1}}}.
        \end{split}
    \end{equation}
     By employing~\eqref{eq:lem-LS-norm:normbigR-sum} to bound \(\normFBig{\bmat{\phi \widehat R_{1:p+1, 1} & \widehat R_{1:p+1, 2:p+1} D_{1:p}^{-1}}}\), we obtain
    \begin{equation} \label{eq:lem-LS:sigminbigR}
        \begin{split}
            & \sigmin\big([\phi \tilde R_{1:p+1, 1} \quad \tilde R_{1:p+1, 2:p+1} D_{1:p}^{-1}]\big) \\
            & \quad \leq (\alpha_{11} + \epsls \alpha_{21}) \normFBig{\bmat{\phi r, \widehat{W}_{1:p} D_{1:p}^{-1}}}
            + (\alpha_{12} + \epsls \alpha_{22}) \normFbig{\widehat{W}_{1:p} D_{1:p}^{-1}} \\
            & \quad \leq (\alpha_{11} + \epsls \alpha_{21}) \phi \beta
            + (\alpha_{11} + \epsls \alpha_{21} + \alpha_{12} + \epsls \alpha_{22}) \normFbig{\widehat{W}_{1:p} D_{1:p}^{-1}},
        \end{split}
    \end{equation}
    which implies that, by multiplying the two sides by \(\phi^{-1}\),
    \begin{equation} \label{eq:lem-LS-sigminbigRDy}
        \begin{split}
            & \sigmin\big([\phi \tilde R_{1:p+1, 1} \quad \tilde R_{1:p+1, 2:p+1} D_{1:p}^{-1}]\big) \phi^{-1} \\
            & \quad \leq (\alpha_{11} + \epsls \alpha_{21}) \beta
            + (\alpha_{11} + \epsls \alpha_{21} + \alpha_{12} + \epsls \alpha_{22}) \normFbig{\widehat{W}_{1:p} D_{1:p}^{-1}} \phi^{-1}.
        \end{split}
    \end{equation}
    Together with~\eqref{eq:lem-LS:norm-Ry}, we have
    \begin{equation} \label{eq:lem-proof:norm-RDR1y}
        \begin{split}
            &\norm{(\widehat{R}_{1:p+1} + \Delta R_{1:p+1}) \bmat{1 \\ -\widehat y^{(k)}}}\\
            &\leq \sqrt{2}\,\bigl((\alpha_{11} + \epsls \alpha_{21}) \beta
            + (\alpha_{11} + \epsls \alpha_{21} + \alpha_{12} + \epsls \alpha_{22}) \normFbig{\widehat{W}_{1:p} D_{1:p}^{-1}} \phi^{-1}\bigr).
        \end{split}
    \end{equation}

    Next, we will prove \(\delta(\phi)\leq 1/2\).
    To bound \(\delta(\phi)\), we will first bound \(\phi\beta\) by \(\normFbig{\widehat{W}_{1:p} D_{1:p}^{-1}}\).
    By~\eqref{eq:sstep-GMRES:beta}, \((\widehat{R}_{1:p+1} + \Delta R_{1:p+1}) \bmat{1 \\ -\widehat y^{(k)}}\) can be written as
    \begin{equation*}
    \begin{split}
        (\widehat{R}_{1:p+1} + &\Delta R_{1:p+1}) \bmat{1 \\ -\widehat y^{(k)}} \\
        ={}& (1 + \epsr) \beta e_1 - \widehat R_{1:p+1, 2:p+1} \widehat y^{(k)} - \Delta R_{1:p+1, 2:p+1} \widehat y^{(k)}.
    \end{split}
    \end{equation*}
    This implies that
    \begin{equation} \label{eq:lem-proof:phibeta-1}
        \begin{split}
            \phi\beta &\leq \frac{\phi}{1 - \epsr} \biggl(\normbigg{(\widehat{R}_{1:p+1} + \Delta R_{1:p+1}) \bmat{1 \\ -\widehat y^{(k)}}} + \normbig{\widehat{R}_{1:p+1, 2:p+1} \widehat y^{(k)}} \\
            &\quad+ \normbig{\Delta R_{1:p+1, 2:p+1} \widehat y^{(k)}}\biggr),
        \end{split}
    \end{equation}
    where \(\normbig{\widehat{R}_{1:p+1, 2:p+1} \widehat y^{(k)}}\) and \(\normbig{\Delta R_{1:p+1, 2:p+1} \widehat y^{(k)}}\) can be bounded by using~\eqref{eq:sstep-GMRES-ls} and Lemma~\ref{lem:LS-norm}:
    \begin{equation} \label{eq:lem-proof:R-DR}
        \begin{split}
            \normbig{\widehat{R}_{1:p+1, 2:p+1} \widehat y^{(k)}}
            &+ \normbig{\Delta R_{1:p+1, 2:p+1} \widehat y^{(k)}} \\
            \leq{}& \normbig{\widehat{R}_{1:p+1, 2:p+1}D_{1:p}^{-1}D_{1:p} \widehat y^{(k)}}
            + \normbig{\Delta R_{1:p+1, 2:p+1} D_{1:p}^{-1}D_{1:p}\widehat y^{(k)}} \\
            \leq{}& (1+\epsls) \normbig{\widehat{R}_{1:p+1, 2:p+1}D_{1:p}^{-1}}\norm{D_{1:p} \widehat y^{(k)}} \\
            \leq{}& (1+\epsls) \alpha_4 \normFbig{\widehat{W}_{1:p} D_{1:p}^{-1}} \norm{D_{1:p} \widehat y^{(k)}}.
        \end{split}
    \end{equation}
    By substituting~\eqref{eq:lem-proof:R-DR} and~\eqref{eq:lem-proof:norm-RDR1y} into~\eqref{eq:lem-proof:phibeta-1}, we obtain
    \begin{equation}
    \begin{split}
        \phi\beta &\leq \frac{1}{1-\epsr}\biggl(\sqrt{2}\,\bigl((\alpha_{11} + \epsls \alpha_{21}) \phi \beta
        + (\alpha_{11} + \epsls \alpha_{21} + \alpha_{12} + \epsls \alpha_{22}) \normFbig{\widehat{W}_{1:p} D_{1:p}^{-1}}\bigr)\\
        &\quad+ (1+\epsls) \alpha_4 \normFbig{\widehat{W}_{1:p} D_{1:p}^{-1}} \norm{D_{1:p} \widehat y^{(k)}}\phi\biggr) \\
        &\leq \frac{1}{1-\epsr}\biggl(\sqrt{2}\,\bigl((\alpha_{11} + \epsls \alpha_{21}) \phi \beta
        + (\alpha_{11} + \epsls \alpha_{21} + \alpha_{12} + \epsls \alpha_{22}) \normFbig{\widehat{W}_{1:p} D_{1:p}^{-1}}\bigr)\\
        &\quad+ (1+\epsls) \alpha_4 \normFbig{\widehat{W}_{1:p} D_{1:p}^{-1}}\biggr).
    \end{split}
    \end{equation}
    The last inequality is derived by noticing \( \norm{D_{1:p} \widehat y^{(k)}}\phi = \sqrt{1-\delta^2(\phi)}\leq 1\) from~\eqref{eq:lem-LS:defphi} and \(0\leq \delta(\phi)<1\). 
    Thus, we bound \(\phi \beta\) by \(\normFbig{\widehat{W}_{1:p} D_{1:p}^{-1}}\) as follows:
    \begin{equation} \label{eq:lem-proof:phibeta}
        \phi \beta \leq \frac{\sqrt{2}\, (\alpha_{11} + \alpha_{12}) + \sqrt{2}\, \epsls (\alpha_{21} + \alpha_{22}) + (1 + \epsls) \alpha_4}
        {1 - \epsr - \sqrt{2}\,\alpha_{11} - \sqrt{2}\, \epsls \alpha_{21}}
        \cdot\normFbig{\widehat{W}_{1:p} D_{1:p}^{-1}}.
    \end{equation}
    Combining~\eqref{eq:lem-proof:phibeta} and~\eqref{eq:lem-LS:sigminbigR} with Lemma~\ref{lem:LS-norm}, we can bound \(\delta(\phi)\) defined by~\eqref{eq:lem-LS:defphi} by
    \begin{equation} \label{eq:lem-LS:deltaphi}
        \begin{split}
            \delta(\phi) \leq{}& \frac{(\alpha_{11} + \epsls \alpha_{21}) \phi \beta
            + (\alpha_{11} + \epsls \alpha_{21} + \alpha_{12} + \epsls \alpha_{22}) \normFbig{\widehat{W}_{1:p} D_{1:p}^{-1}}}
            {(1-\epsls)\alpha_3 \normFbig{\widehat{W}_{1:p} D_{1:p}^{-1}}} \\
            \leq{}& \frac{(\alpha_{11} + \epsls \alpha_{21}) \frac{\sqrt{2}\, (\alpha_{11} + \alpha_{12}) + \sqrt{2}\, \epsls (\alpha_{21} + \alpha_{22}) + \alpha_4}
            {1 - \epsr - \sqrt{2}\, \alpha_{11} - \sqrt{2}\, \epsls \alpha_{21}}
            + (\alpha_{11} + \epsls \alpha_{21} + \alpha_{12} + \epsls \alpha_{22})}
            {(1-\epsls)\alpha_3}.
        \end{split}
    \end{equation}
    Utilizing the assumption~\eqref{eq:lem-LS:assump} and ignoring the quadratic terms, it can be checked that \(\delta(\phi)\leq 1/2\).
    
    Then using~\eqref{eq:lem-proof:norm-RDR1y} along with~\(\delta(\phi)\leq 1/2\) and~\eqref{eq:lem-proof:R-DR}, it follows that
    \begin{equation*}
        \begin{split}
            \normbigg{\widehat{R}_{1:p+1} \bmat{1 \\ -\widehat y^{(k)}}}
            \leq{}& \normbigg{(\widehat{R}_{1:p+1} + \Delta R_{1:p+1}) \bmat{1 \\ -\widehat y^{(k)}}} + \normbigg{\Delta R_{1:p+1} \bmat{1 \\ -\widehat y^{(k)}}} \\
            \leq{}& (3 \cdot \alpha_{11} + 4 \cdot \epsls \alpha_{21}) \beta
            + (3 \cdot \alpha_{11} + 4 \cdot \epsls \alpha_{21} + 3 \cdot \alpha_{12} \\
            & + 4 \cdot \epsls \alpha_{22}) \normFbig{\widehat{W}_{1:p} D_{1:p}^{-1}} \norm{D_{1:p} \widehat{y}^{(k)}},
        \end{split}
    \end{equation*}
    which concludes the proof by substituting \(\alpha_*\) from Lemma~\ref{lem:LS-norm} and ignoring the quadratic terms.
\end{proof}

Lemma~\ref{lem:LS} not only bounds the residual of the least squares problem, but it also establishes the connection between \(\epskappaW\) and \(\epskapparW\), indicating that only \(\epskapparW\) needs to be determined.

The following theorem gives the upper bound on \(\norm{b - A \widehat x^{(k)}}\) and illustrates how the backward error is affected by the errors from various steps.

\begin{theorem} \label{thm:b-Axi-H}
    Assume that \(\widehat{W}_{1:is}\), \(\widehat{V}_{is+1}\), \(\widehat{R}_{is+1, is+1}\), \(\widehat{G}_{1:is, 1:is}\), \(\widehat{T}_{1:is, 1:is}\), \(\widehat{y}^{(i)}\), and \(\widehat{x}^{(i)}\) satisfy~\eqref{eq:sstep-GMRES-wi}--\eqref{eq:sstep-GMRES-VR}, \eqref{eq:sstep-GMRES-ls}, and \eqref{eq:sstep-GMRES-xi}.
    If there exists an iteration \(p := \ki s + \kj \leq n\) such that~\eqref{eq:lem:LS:assump-rW} and
    \begin{equation} \label{eq:thm-b-Axi-H:assump}
        \epskappaW \frac{\sqrt{p}\, (1 + \epsAv + \epsMLAZ) \kappa(M_L) \kappa(A) \kappa(\tilde Z_{1:p})}{1 - \sqrt{p}\, (\epsAv + \epsMLAZ) \kappa(M_L) \kappa(A) \kappa(\tilde Z_{1:p})} \leq 1
    \end{equation}
    hold with~\eqref{eq:lem-LS:assump},
    then
    \begin{equation} \label{eq:lem:b-Axi-1}
        \begin{split}
            \norm{b - A\widehat{x}^{(\ki)}}
            &\leq \xi_b \kappa(M_L) \norm{b} + \xi_{Ax} \frac{\sqrt{p}\, \kappa (\tilde Z_{1:p})}{1 - \sqrt{p}\, \epsMV \kappa (\tilde Z_{1:p})} \kappa(M_L) \normF{A} \norm{\widehat x^{(\ki)}}
        \end{split}
    \end{equation}
    with
    \begin{align*}
        \xi_b ={}& \sqrt{1 + \omega_{p+1}}\, (9 \cdot \epskapparW + 10 \cdot \epsqr + 12 \cdot \epsls) (1 + 2 \cdot \epsMLb), \\
        \xi_{Ax} ={}& \sqrt{1 + \omega_{p+1}}\, (1 + \epsAv + \epsMLAZ)
        (9 \cdot \epskapparW + 10 \cdot \epsqr + 12 \cdot \epsls + 12 \cdot \epsH) \\
        & + \epsAv + \epsMLAZ + \epsMV,
    \end{align*}
    where \(\widehat Z_{1:p} = \tilde Z_{1:p} D_{1:p}\) with a positive definite diagonal matrix \(D_{1:p}\), and
    \[
        \kappa(\tilde Z_{1:p}) \leq \frac{(1 + \sqrt{p}\, \epsB) \kappa(M_R) \kappa(\tilde B_{1:p})}{1 - \sqrt{p}\, \epsB \kappa(M_R) \kappa(\tilde B_{1:p})}
    \]
    with \(\widehat B_{1:p} = \tilde B_{1:p} D_{1:p}\).
\end{theorem}

\begin{proof}
    For the \(\ki\)-th block iteration, we will estimate the preconditioned residual \(\norm{M_L^{-1} (b - A \widehat{x}^{(\ki)})}\), from which we can easily derive the bound of the residual \(\norm{b - A\widehat{x}^{(\ki)}}\).
    From~\eqref{eq:sstep-GMRES-wi} and \eqref{eq:sstep-GMRES-VR}, we summarize the error from generating the basis and performing the block orthogonalization as
    \begin{equation} \label{eq:lem:b-Axi-H:rML-1AZ}
        \begin{split}
        \bmat{\widehat{r} & M_L^{-1} A \widehat{Z}_{1:p}}
        &= \bmat{\widehat{r} & \widehat{W}_{1:p}} - \bmat{0& \Delta W_{1:p}} \\
        &= \bar{V}_{1:p+1} \widehat{R}_{1:p+1} (1) - \underbrace{\bigl(\Delta E_{1:p+1} (1) + \bmat{0& \Delta W_{1:p}}\bigr)}_{=:\Delta F_{1:p+1}},
        \end{split}
    \end{equation}
    where \(\Delta F_{1:p+1}\) satisfies, from~\eqref{eq:sstep-GMRES-wi} and~\eqref{eq:sstep-GMRES-DEyi},
    \begin{equation} \label{eq:lem:b-Axi-H:DF}
        \begin{split}
            \normbigg{\Delta F_{1:p+1}\bmat{1\\ -\widehat{y}^{(\ki)}}} \leq{}& \normbigg{\Delta E_{1:p+1} (1) \bmat{1\\ -\widehat{y}^{(\ki)}}} + \normbigg{\bmat{0& \Delta W_{1:p}}\bmat{1\\ -\widehat{y}^{(\ki)}}} \\
            \leq{}& \epsqr \beta + \bigl(\epsqr \bigl(1 + \epsAv + \epsMLAZ\bigr) + \epsAv + \epsMLAZ\bigr) \\
            & \cdot  \norm{M_L^{-1}} \normF{A} \normF{\tilde Z_{1:p}} \norm{D_{1:p} \widehat{y}^{(\ki)}}.
        \end{split}
    \end{equation}
    Using \eqref{eq:sstep-GMRES-xi} along with~\eqref{eq:sstep-GMRES:barr}, \(\norm{M_L^{-1} (b - A \widehat{x}^{(\ki)})}\) can be bounded by
    \begin{equation*} \label{eq:lem:b-Axi-0}
        \begin{split}
            & \norm{M_L^{-1} (b - A \widehat{x}^{(\ki)})} \\
            &\quad \leq \norm{M_L^{-1} b - M_L^{-1} A \widehat{Z}_{1:p} \widehat{y}^{(\ki)}} + \norm{M_L^{-1} A \Delta x^{(\ki)}} \\
            &\quad \leq \normbigg{\bmat{\widehat{r} & M_L^{-1} A \widehat{Z}_{1:p}} \bmat{1\\ -\widehat{y}^{(\ki)}}} + \norm{\Delta r} + \norm{M_L^{-1}} \norm{A} \norm{\Delta x^{(\ki)}}.
        \end{split}
    \end{equation*}
    Furthermore, utilizing the bound of \(\normbigg{\bmat{\widehat{r} & M_L^{-1} A \widehat{Z}_{1:p}} \bmat{1\\ -\widehat{y}^{(\ki)}}}\) derived from~\eqref{eq:lem:b-Axi-H:rML-1AZ} and the bound of \(\bar{V}_{1:p+1}\) from~\eqref{eq:sstep-GMRES-Viorth}, we have
    \begin{equation*}
        \begin{split}
            & \norm{M_L^{-1} (b - A \widehat{x}^{(\ki)})} \\
            &\quad \leq \normbigg{\bar{V}_{1:p+1} \widehat{R}_{1:p+1} (1) \bmat{1\\ -\widehat{y}^{(\ki)}}} + \normbigg{\Delta F_{1:p+1}\bmat{1\\ -\widehat{y}^{(\ki)}}} + \norm{\Delta r} \\
            &\qquad + \norm{M_L^{-1}} \norm{A} \norm{\Delta x^{(\ki)}} \\
            &\quad \leq \sqrt{1 + \omega_{p+1}}\, \normbigg{\widehat{R}_{1:p+1} (1) \bmat{1\\ -\widehat{y}^{(\ki)}}} + \normbigg{\Delta F_{1:p+1}\bmat{1\\ -\widehat{y}^{(\ki)}}} \\
            &\qquad + \epsMLb \norm{M_L^{-1}} \norm{b} + \epsMV \norm{M_L^{-1}} \normF{A} \normF{\tilde Z_{1:p}} \norm{D_{1:p} \widehat{y}^{(\ki)}},
        \end{split}
    \end{equation*}
    which further implies, by~\eqref{eq:lem:b-Axi-H:DF} and Lemma~\ref{lem:LS},
    \begin{equation} \label{eq:lem:b-Axi}
        \begin{split}
            & \norm{M_L^{-1} (b - A \widehat{x}^{(\ki)})} \\
            &\quad \leq \sqrt{1 + \omega_{p+1}}\, (9 \cdot \epskapparW + 10 \cdot \epsqr + 12 \cdot \epsls) (1 + 2 \cdot \epsMLb) \norm{M_L^{-1}} \norm{b} \\
            &\qquad + \bigl(\sqrt{1 + \omega_{p+1}}\, (1 + \epsAv + \epsMLAZ) \cdot (9 \cdot \epskapparW + 10 \cdot \epsqr + 12 \cdot \epsls \\
            &\qquad + 12 \cdot \epsH) + \epsAv + \epsMLAZ + \epsMV \bigr) \norm{M_L^{-1}} \normF{A} \normF{\tilde Z_{1:p}} \norm{D_{1:p} \widehat{y}^{(\ki)}}.
        \end{split}
    \end{equation}
    Note that using Lemma~\ref{lem:LS} requires~\eqref{eq:lem:LS:assump-W} which is guaranteed by the assumption~\eqref{eq:thm-b-Axi-H:assump}.

    Then it remains to bound \(\norm{D_{1:p} \widehat{y}^{(\ki)}}\).
    By
    \begin{equation*}
        \begin{split}
            \sigmin(\tilde Z_{1:p}) \norm{D_{1:p} \widehat y^{(\ki)}} & \leq \norm{\tilde Z_{1:p} (D_{1:p} \widehat y^{(\ki)})} \\
            & \leq \norm{x^{(\ki)}} + \norm{\Delta x^{(\ki)}} \\
            & \leq \norm{x^{(\ki)}} + \epsMV \normF{\tilde Z_{1:p}} \norm{D_{1:p} \widehat y^{(\ki)})},
        \end{split}
    \end{equation*}
    we then derive
    \begin{equation*} 
        \norm{D_{1:p} \widehat y^{(\ki)}} \leq \frac{\norm{\widehat x^{(\ki)}}}{\sigmin(\tilde Z_{1:p}) - \epsMV \normF{\tilde Z_{1:p}}} 
        = \frac{\norm{\widehat x^{(\ki)}}}{\sigmin(\tilde Z_{1:p}) \bigl(1 - \sqrt{p}\, \epsMV \kappa (\tilde Z_{1:p})\bigr)},
    \end{equation*}
    which proves the conclusion combined with~\eqref{eq:lem:b-Axi} and \(\norm{b - A \widehat{x}^{(\ki)}} \leq \norm{M_L} \norm{M_L^{-1} (b - A \widehat{x}^{(\ki)})}\).
\end{proof}

For simplicity, Theorem~\ref{thm:b-Axi-H} does not consider the restarted \(s\)-step GMRES algorithm.
Analogous to~\cite[Theorem 4.1]{BHMV2024}, it is easy to generalize Theorem~\ref{thm:b-Axi-H} to the restarted \(s\)-step GMRES algorithm, which can be regarded as an iterative refinement to remove \(\kappa(M_L)\), \(\kappa(M_R)\), and \(\kappa(\tilde B_{1:p})\) from the bound of the backward error in Theorem~\ref{thm:b-Axi-H} under certain conditions.

\begin{remark} \label{remark:thm-backward}
According to Remark~\ref{remark:eps-1}, \(\epsr\), \(\epsAv\), \(\epsls\), \(\epsMV = \bigO(\macheps)\), and the \(\delta_*\) terms are highly preconditioner dependent.
It remains to evaluate \(1 + \omega_{p+1}\), \(\epsqr\), \(\epskappaW\), \(\epskapparW\), and \(\epsH\), which are all related only to the block orthogonalization method.
Among these terms, it is usually easy to check that \(1 + \omega_{p+1} \leq \bigO(\sqrt{n}\,)\) and \(\epsqr = \bigO(\macheps)\) for different block orthogonalization methods.
In addition, \(\epskappaW\) can be obtained from the connection with \(\epskapparW\) as shown in Lemma~\ref{lem:LS}.
Determining \(\epskapparW\) is crucial for the backward error, since Theorem~\ref{thm:b-Axi-H} demonstrates that \(\epskapparW\) directly affects the backward error, while \(\epskappaW\) determines the possibility of achieving the backward error specified in equation (3.38) for the given linear system.

Then from Theorem~\ref{thm:b-Axi-H} the difficulty of analyzing the backward error of the \(s\)-step GMRES algorithm is to prove that~\eqref{eq:lem:LS:assump-rW} is satisfied when \(\widehat V_{1:p+1}\) is not well-conditioned, namely, to determine \(\epskapparW\).
In the case that \(\bar V_*\) is not exactly orthonormal, we need to additionally estimate \(\epsH\).
\end{remark}

We now turn to the analysis of the backward stability of $s$-step GMRES with specific block orthogonalization schemes. As we will see for all orthogonalization methods in Section~\ref{sec:ortho}, there is a dependence of the backward error on \(\kappa(\tilde{B}_{i:p})\). It is not clear how to bound this quantity a priori.
We will give an example in Section~\ref{sec:modifiedArnoldi} that shows that \(\kappa(\tilde{B}_{i:p})\) can be large even when the condition numbers of its blocks are small, and then we will present a modified approach to eliminate this problem.

\section{Backward stability of \(s\)-step GMRES with different orthogonalization methods}
\label{sec:ortho}
Based on the framework introduced in Section~\ref{sec:stability}, we now analyze the \(s\)-step GMRES algorithm combined with specific block orthogonalization methods.
Note that for simplicity, we assume \(M_L = I\) and \(M_R = I\) in this section, because the constants related to preconditioners vary significantly depending on the particular preconditioners chosen.
Consequently, \(\kappa(M_L) = \kappa(M_R) = 1\), \(\kappa(\tilde{Z}) = \kappa(\tilde{B})\), and \(\epsMLb = \epsMLAZ = 0\).
For readers who wish to analyze the backward error when using specific preconditioners, only \(\kappa(M_L)\), \(\kappa(\tilde{Z})\), \(\epsMLb\), and \(\epsMLAZ\) in Theorem~\ref{thm:b-Axi-H} must be determined.

The bounds for the relative backward error derived in this section will all ultimately depend on the condition number of the computed basis for the Krylov subspace. This confirms what has been widely observed experimentally, i.e., that the conditioning of the Krylov basis is crucial in determining the resulting numerical behavior of $s$-step GMRES. It is, unfortunately, difficult to bound this quantity a priori. In Section \ref{sec:modifiedArnoldi}, we comment on attempting to control the condition number and present a modified approach for ensuring that this quantity remains small. 

\subsection{Backward error of \(s\)-step GMRES with block Householder QR and block modified Gram--Schmidt (BMGS)}
\label{subsec:GMRES-HH-bmgs}
In~\cite{Walker1988}, Walker discussed the GMRES algorithm with Householder QR orthogonalization, which can be directly extended to \(s\)-step GMRES with block Householder QR.
Note that using Householder QR in GMRES is more computationally intensive than alternative orthogonalization schemes.
As a result, there is limited research on employing block Householder QR in the \(s\)-step GMRES algorithm.
However, for purposes of illustration, we also analyze the backward stability of this variant. 

From~\cite[Theorem 19.4]{H2002} and~\cite{YFS2021}, it is easy to see that
\begin{equation} \label{eq:sstep-GMRES-HH-VR}
    \bmat{\phi \widehat{r}, \widehat{W}_{1:n-1}} + \Delta E_{1:n}(\phi) = \tilde V_{1:n} \widehat{R}_{1:n}(\phi), \quad \norm{\Delta E_j(\phi)} \leq \bigO(\macheps) \normBig{\bmat{\widehat{r}, \widehat{W}_{1:n-1}}_j},
\end{equation}
for any \(j \in \{1, \dotsc, n\}\),
is satisfied with an exact orthogonal matrix \(\tilde V_{1:n}\) for the block Householder QR algorithm.
This means that for block Householder QR orthogonalization, \eqref{eq:lem:LS:assump-rW} is satisfied with \(\epskapparW = \bigO(\macheps)\)  when \(p = n\), since the block Householder QR algorithm generates a well-conditioned set of vectors until \(p = n\).

Then we consider the modified Gram--Schmidt (MGS) algorithm.
\cite[Theorem 19.13]{H2002} and~\cite{BP1992} showed~\eqref{eq:sstep-GMRES-HH-VR} can be satisfied by establishing the equivalence between the MGS algorithm and the Householder QR algorithm.
Similarly, it has been proved in~\cite{B2019} that applying the BMGS algorithm with  Householder QR as the intra-block orthogonalization routine to \(X\) is equivalent to applying the block Householder QR algorithm to \(\bmat{0_n \\ X}\) both mathematically and numerically.
Then we obtain that~\eqref{eq:sstep-GMRES-HH-VR} holds also for the BMGS algorithm.

Unlike the block Householder QR algorithm, we cannot ensure the orthogonality of \(\widehat V_{1:n}\).
Therefore, it is necessary to demonstrate that \(\normF{\widehat V_{1:j}\trans \widehat V_{1:j} - I} \leq 1/2\) if 
\[
\sigmin\Big(\bmat{\phi \widehat{r}, \widehat{W}_{1:j} D_{1:j}^{-1}}\Big) / 2 \geq \bigO(\macheps) \normFBig{\bmat{\phi \widehat{r}, \widehat{W}_{1:j} D_{1:j}^{-1}}}
\] 
holds, as established in~\cite[Theorem~4.1]{JP1991}.
Moreover, the contrapositive indicates that~\eqref{eq:lem:LS:assump-rW} is satisfied with \(\epskapparW = \bigO(\macheps)\) when there exists \(p \leq n\) such that \(\normF{\widehat V_{1:p+1}\trans \widehat V_{1:p+1} - I} > 1/2\) but \(\normF{\widehat V_{1:p}\trans \widehat V_{1:p} - I} \leq 1/2\).

Therefore, combined with \(\epsr\), \(\epsqr\), \(\epsAv\), \(\epsls\), \(\epsMV = \bigO(\macheps)\), and \(\omega_p = 0\) from Remark~\ref{remark:eps-1}, Theorem~\ref{thm:b-Axi-H} implies the following lemma.

\begin{lemma} \label{lem:sstep-GMERS-hh-bmgs}
    Assume that \(\widehat{W}_{1:is}\), \(\widehat{V}_{is+1}\), \(\widehat{R}_{is+1, is+1}\), \(\widehat{G}_{1:is, 1:is}\), \(\widehat{T}_{1:is, 1:is}\), \(\widehat{y}^{(i)}\), and \(\widehat{x}^{(i)}\) satisfying~\eqref{eq:sstep-GMRES-wi}--\eqref{eq:sstep-GMRES-VR}, \eqref{eq:sstep-GMRES-ls}, and \eqref{eq:sstep-GMRES-xi}, are computed by the \(s\)-step GMRES algorithm with block Householder QR or BMGS orthogonalization.
    There exists \(p = \ki s + \kj\) such that~\eqref{eq:lem:LS:assump-rW} holds.
    If it also holds that
    \begin{equation*}
        \frac{\bigO(\macheps) \kappa(A) \kappa(\tilde B_{1:p})}{1 - \bigO(\macheps) \kappa(A) \kappa(\tilde B_{1:p})} \leq 1,
    \end{equation*}
    then \(\kappa(\widehat V_{1:p}) \leq 3\) and
    \begin{equation*}
        \begin{split}
            \frac{\norm{b - A\widehat{x}^{(\ki)}}}{\norm{b} + \normF{A} \norm{\widehat x^{(\ki)}}}
            \leq{}& \frac{\bigO(\macheps) \kappa (\tilde B_{1:p})}{1 - \bigO(\macheps) \kappa (\tilde B_{1:p})},
        \end{split}
    \end{equation*}
    where \(\widehat B_{1:p} = \tilde B_{1:p} D_{1:p}\) with a positive definite diagonal matrix \(D_{1:p}\).

    Furthermore, if \(s = 1\), as long as
    \(
        \frac{\bigO(\macheps) \kappa(A)}{1 - \bigO(\macheps) \kappa(A)} \leq 1,
    \)
    then there exists \(p = \ki s + \kj\) such that \(\kappa (\widehat B_{1:p}) = \kappa (\widehat V_{1:p}) \leq 3\) and
    \begin{equation*}
        \begin{split}
            \frac{\norm{b - A\widehat{x}^{(\ki)}}}{\norm{b} + \normF{A} \norm{\widehat x^{(\ki)}}}
            \leq{}& \frac{\bigO(\macheps)}{1 - \bigO(\macheps)}.
        \end{split}
    \end{equation*}
\end{lemma}

For \(s = 1\), i.e., the standard GMRES algorithm with Householder QR and MGS orthogonalization, we recover the same backward stability result of~\cite{DGRS1995} and~\cite{PRS2006}, respectively.

\subsection{Backward error of \(s\)-step GMRES with  reorthogonalized block classical Gram--Schmidt (BCGSI+)}
In this subsection, we examine another scenario, such as BCGSI+.
As derived from~\cite[Corollary 1]{CM2024-stableonesync}, it is evident that we can achieve~\eqref{eq:sstep-GMRES-HH-VR} using an exactly orthonormal matrix \(\tilde V\) for BCGSI+.
Nevertheless, in this case, we also offer a more complicated proof under the assumption that obtaining~\eqref{eq:sstep-GMRES-HH-VR} with an exactly orthonormal matrix \(\tilde V\) is not feasible, as detailed in Appendix~\ref{appendix:bcgsi+}.
This analysis may shed light on orthogonalization methods that indeed cannot accomplish~\eqref{eq:sstep-GMRES-HH-VR} with an exactly orthonormal matrix \(\tilde V\).
According to Lemma~\ref{lem:LS}, we need to prove both~\eqref{eq:lem:LS:assump-rW} and~\eqref{eq:lem-LS:assump-epsH}.
We provide the properties of BCGSI+ in Appendix~\ref{appendix:bcgsi+}.

By Lemma~\ref{lem:bcgs2:backwarderror-loo}, we obtain that \(\epsqr = \bigO(\macheps)\) and~\eqref{eq:lem:LS:assump-rW} holds with \(\epskapparW = \bigO(\macheps)\).
From Lemma~\ref{lem:epsH} and the fact that \(\widehat T\) is upper triangular, \eqref{eq:lem-LS:assump-epsH} holds with \(\epsH = \bigO(\macheps)\).
Together with \(\epsr\), \(\epsAv\), \(\epsls\), \(\epsMV = \bigO(\macheps)\) as discussed in Remark~\ref{remark:eps-1}, we derive the following lemma to show the backward stability of the \(s\)-step GMRES with the BCGSI+ algorithm using Theorem~\ref{thm:b-Axi-H}.

\begin{lemma} \label{lem:sstep-GMERS-bcgs2}
    Assume that \(\widehat{W}_{1:is}\), \(\widehat{V}_{is+1}\), \(\widehat{R}_{is+1, is+1}\), \(\widehat{G}_{1:is, 1:is}\), \(\widehat{T}_{1:is, 1:is}\), \(\widehat{y}^{(i)}\), and \(\widehat{x}^{(i)}\) satisfying~\eqref{eq:sstep-GMRES-wi}--\eqref{eq:sstep-GMRES-VR}, \eqref{eq:sstep-GMRES-ls}, and \eqref{eq:sstep-GMRES-xi}, are computed by the \(s\)-step GMRES algorithm with BCGSI+ orthogonalization.
    There exists \(p = \ki s + \kj\) such that~\eqref{eq:lem:LS:assump-rW} holds.
    If it also holds that
    \begin{equation*}
        \frac{\bigO(\macheps) \kappa(A) \kappa(\tilde B_{1:p})}{1 - \bigO(\macheps) \kappa(A) \kappa(\tilde B_{1:p})} \leq 1,
    \end{equation*}
    then \(\kappa(\widehat V_{1:p}) \leq \frac{1 + \bigO(\macheps)}{1 - \bigO(\macheps)}\) and
    \begin{equation*}
        \begin{split}
            \frac{\norm{b - A\widehat{x}^{(\ki)}}}{\norm{b} + \normF{A} \norm{\widehat x^{(\ki)}}}
            \leq{}& \frac{\bigO(\macheps) \kappa (\tilde B_{1:p})}{1 - \bigO(\macheps) \kappa (\tilde B_{1:p})},
        \end{split}
    \end{equation*}    
    where \(\widehat B_{1:p} = \tilde B_{1:p} D_{1:p}\) with a positive definite diagonal matrix \(D_{1:p}\).

    Furthermore, if \(s = 1\), as long as
    \(
        \frac{\bigO(\macheps) \kappa(A)}{1 - \bigO(\macheps) \kappa(A)} \leq 1,
    \)
    then there exists \(p = \ki s + \kj\) such that \(\kappa (\widehat B_{1:p}) = \kappa (\widehat V_{1:p}) \leq \frac{1 + \bigO(\macheps)}{1 - \bigO(\macheps)}\) and
    \begin{equation*}
        \begin{split}
            \frac{\norm{b - A\widehat{x}^{(\ki)}}}{\norm{b} + \normF{A} \norm{\widehat x^{(\ki)}}}
            \leq{}& \frac{\bigO(\macheps)}{1 - \bigO(\macheps)}.
        \end{split}
    \end{equation*}
\end{lemma}

For \(s = 1\), i.e., the standard GMRES algorithm with CGSI+ orthogonalization, we recover the same backward stability result of \cite{DGRS1995}.

\section{Discussion of theoretical results}
\label{sec:discussion}
In this part, we discuss the stopping criteria and the requirements of the orthogonalization method indicated by the above theoretical results.

\subsection{Stopping criteria}
The commonly used stopping criteria for the GMRES algorithm are
\begin{equation} \label{eq:stopping-ls}
    \widehat \beta \widehat G_{1, is+1} \leq \texttt{tolLS} \cdot \norm{r}
\end{equation}
and
\begin{equation} \label{eq:stopping-backward-err}
    \norm{b - A \widehat{x}^{(i)}} \leq \texttt{tol} \cdot \bigl(\norm{b} + \normF{A} \norm{\widehat x^{(i)}}\bigr),
\end{equation}
where \texttt{tolLS} and \texttt{tol} are user-specified thresholds. 
The first criterion~\eqref{eq:stopping-ls} is straightforward and economical to check practically, yet it might fail to recognize timely convergence.
This is because \(\widehat \beta \widehat G_{1, is+1}/\norm{r}\) might not be sufficiently small, even when \(\norm{b - A \widehat{x}^{(i)}}/\bigl(\norm{b} + \normF{A} \norm{\widehat x^{(i)}}\bigr)\) is sufficiently small.
The disadvantage of~\eqref{eq:stopping-backward-err} is that it requires more computational effort, as \(\widehat{x}^{(i)}\) and \(A \widehat{x}^{(i)}\) must be determined.
Consequently, in practice, both criteria are generally employed.
This implies that ~\eqref{eq:stopping-ls} is checked in every iteration, whereas ~\eqref{eq:stopping-backward-err} is evaluated in select iterations when implementing the GMRES algorithm.

For the \(s\)-step GMRES with \(s = 1\), our theoretical results, Lemmas~\ref{lem:sstep-GMERS-hh-bmgs} and~\ref{lem:sstep-GMERS-bcgs2}, show that the backward error is \(\bigO(\macheps)\), provided that \(A\) is numerically nonsingular.
This aligns with the existing results for standard GMRES.
However, for \(s > 1\), the situation differs.
From Theorem~\ref{thm:b-Axi-H}, as well as Lemmas~\ref{lem:sstep-GMERS-hh-bmgs} and~\ref{lem:sstep-GMERS-bcgs2}, we find that there is no guarantee for the backward error of the \(s\)-step GMRES algorithm with \(s > 1\) due to the existence of \(\kappa (\tilde B_{1:p})\).
This implies that for the ``key dimension'' \(p\), the \(s\)-step GMRES algorithm nearly reaches optimal accuracy when~\eqref{eq:lem:LS:assump-rW} is satisfied.
Meanwhile, if~\eqref{eq:lem:LS:assump-W} also holds, the backward error can be bounded by \(\bigO(\macheps) \kappa (\tilde B_{1:p})\).
Otherwise, there is no theoretical guarantee for this ``optimal'' accuracy.
This means that it is not possible to find this ``optimal'' accuracy only through detecting~\eqref{eq:stopping-ls} and~\eqref{eq:stopping-backward-err}.

Based on our analysis, the time when \(\widehat V_{1:p+1}\) loses orthogonality also indicates that
\[
    \frac{\max\bigl(\abs{\widehat R_{p+1, p+1} D_{p, p}^{-1}}, \normF{\bar V_{p+1} \widehat R_{p+1, 1:p+1} D_{1:p}^{-1}}\bigr)}{\normF{\widehat W_{1:p} D_{1:p}^{-1}}}
\]
is sufficiently small, as described in~\eqref{eq:lem-LS:assump-epsH}.
Therefore, we introduce an additional criterion to identify the ``key dimension'' by
\begin{equation} \label{eq:stoppingg-H}
    \abs{\widehat R_{p+1, p+1}} \leq \texttt{tolH} \cdot \normF{\widehat W_{1:p}},
\end{equation}
for a user-specified threshold \texttt{tolH}, 
which can terminate the algorithm when it achieves approximately ``optimal'' accuracy.
In standard GMRES, \eqref{eq:stoppingg-H} is rarely used, since \eqref{eq:stopping-backward-err} or~\eqref{eq:stopping-ls} is satisfied but \eqref{eq:stoppingg-H} is not for many cases.
The reason for this is that \eqref{eq:stoppingg-H} indicates the quality of the approximation of \(A^{-1}\), which is a challenge to estimate.

\subsection{Requirement of the orthogonalization method}
In Theorem~\ref{thm:b-Axi-H} and Remark~\ref{remark:thm-backward}, we show that \(\epskapparW\) is determined by the orthogonalization method and has a direct effect on the backward error~\eqref{eq:lem:b-Axi-1}.
Revisiting the above proof, we bound \(\epskapparW\) by analyzing the loss of orthogonality of the orthogonalization method used for the QR factorization of \(X\), specifically \(X = QR\).
If the orthogonalization method satisfies 
\begin{equation}
    \normF{\widehat Q\trans \widehat Q - I} \leq \frac{1}{2},
\end{equation}
provided \(\bigO(\macheps) \kappa^{\alpha} (X) \leq 1\), then the contrapositive
indicates that \(\bigO(\macheps) \kappa^{\alpha} (X) > 1\) when \(\widehat Q\) is not well-conditioned.
Note that when \(\widehat{Q}\) is not well-conditioned, it implies that GMRES reaches the key dimension \(p\).
As mentioned in Section~\ref{sec:stability}, \(\epskapparW\), as defined in~\eqref{eq:lem:LS:assump-rW}, can be deduced from \(\bigO(\macheps) \kappa^{\alpha} (X) > 1\), which amounts to \(\sigmin(X) < \bigO(\macheps^{1/\alpha}) \normF{X}\).
By substituting \(\bmat{\widehat{r} & \widehat{W}_{1:p}}\) for \(X\),
it follows that \(\epskapparW = \bigO(\macheps^{1/\alpha})\) in~\eqref{eq:lem:LS:assump-rW}.

Furthermore, even for standard GMRES, i.e., \(s = 1\), employing an orthogonalization method with \(\alpha = 2\), such as the reorthogonalized Pythagorean variants of BCGS introduced by~\cite{CLMO2024-P, CLRT2022}, implies that \(\epskapparW = \bigO(\sqrt{\macheps}\,)\).
Note that other terms, \(\omega_*\) and \(\varepsilon_*\), are usually \(\bigO(\macheps)\), except for \(\epsH\) which also depends on the orthogonalization method.
As a result, ignoring the effect of preconditioning, \(\xi_b\) and \(\xi_{Ax}\) from~\eqref{eq:lem:b-Axi-1} are dominated by \(\epskapparW\), meaning that \(\xi_b\) and \(\xi_{Ax}\) are at least \(\bigO(\sqrt{\macheps}\,)\).
Thus, plugging them into~\eqref{eq:lem:b-Axi-1}, we can only expect the backward error to be bounded by
\begin{equation*}
    \begin{split}
        \frac{\norm{b - A\widehat{x}^{(\ki)}}}{\norm{b} + \normF{A} \norm{\widehat x^{(\ki)}}}
        \leq{}& \bigO(\sqrt{\macheps}\,),
    \end{split}
\end{equation*}
which is illustrated by an example in Appendix~\ref{appendix:example-pip}.

This suggests that we cannot prove a \(\bigO(\macheps)\)-level backward error for GMRES or $s$-step GMRES with  orthogonalization methods that require \(\bigO(\macheps) \kappa^\alpha(X)\) with \(\alpha > 1\).

\section{A modified Arnoldi process for improving the stability of \(s\)-step GMRES}
\label{sec:modifiedArnoldi}
As indicated in Section~\ref{sec:stability}, the backward stability of \(s\)-step GMRES algorithm is directly influenced by the condition number of
\begin{equation*}
    \begin{split}
        \tilde B_{1:p} & = \hat B_{1:p} D^{-1}_{1:p} \\
        & = \bmat{\hat B_{1:s} D^{-1}_{1:s}& \dotsb & \hat B_{(\ki-1)s+1:\ki s} D^{-1}_{(\ki-1)s+1:\ki s}& \hat B_{\ki s+1:p} D^{-1}_{\ki s+1:p}} \\
        & = \bmat{\tilde B_{1:s}& \dotsb & \tilde B_{(\ki-1)s+1:\ki s}& \tilde B_{\ki s+1:p}}.
    \end{split}
\end{equation*}
In the classical \(s\)-step GMRES algorithm, the Krylov submatrices \(K_{(i-1)s+1:is}\) are utilized as \(B_{(i-1)s+1:is}\), with \(K_{(i-1)s+1:is}\) being formed as 
\begin{equation} \label{eq:def-Zi}
K_{(i-1)s+1:is} = \bmat{p_0(A) V_{(i-1)s+1} & p_1(A) V_{(i-1)s+1} & \dotsi & p_{s-1}(A) V_{(i-1)s+1}},
\end{equation}
incorporating the \(s\)-step basis polynomials \(p_0\), \(p_1\), \(\dotsc\), \(p_{s-1}\).
Here, popular choices of polynomials include monomial, Newton~\cite{BHR1994}, and Chebyshev polynomials~\cite{JC1992-1, JC1992-2, DV1995}.

One method to manage \(\kappa (\tilde B_{1:p})\) is to regulate each sub-block of \(\tilde B_{1:p}\) by adaptively selecting different \(s\) in the algorithm, as advocated in~\cite{IE2017}.
This approach is effective in many scenarios.
However, in certain specific cases, \(\kappa (\tilde B_{1:p})\) can be very large even if the condition number of each sub-block of \(\tilde B_{1:p}\) is small.
This implies that merely constraining the condition number of each sub-block is insufficient to bound \(\kappa (\tilde B_{1:p})\). We demonstrate this through an example below. 

\begin{figure}[!tb]
    \centering
    \includegraphics[width=1.0\textwidth]{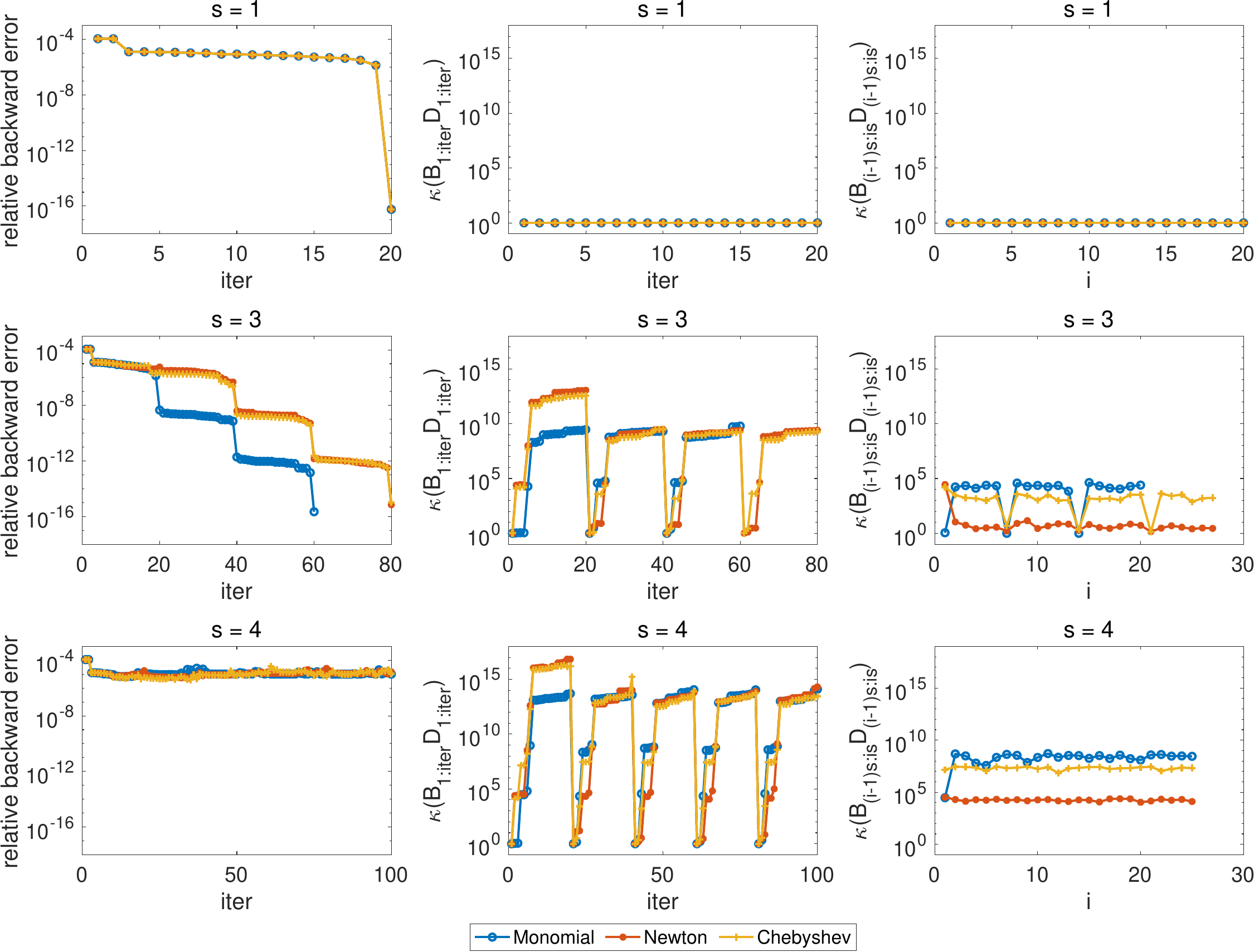}
    \caption{The plot for Example~\ref{example:badcase}: From left to right, the plots are of relative backward error, the condition number of the basis \(\tilde{B}\), and the condition number of the sub-block for the basis of each iteration, where we normalize each column of \(\hat B\) as \(\tilde B\).
    Each line with different color in the plots denotes \(s\)-step GMRES using different polynomials, including monomial, Newton, and Chebyshev polynomials, to generate the basis by~\eqref{eq:def-Zi}.}
    \label{fig:badcase}
\end{figure}

\begin{example} \label{example:badcase}
    We construct the linear system \(Ax = b\), where \(A\) is a \(20\)-by-\(20\) random matrix with \(\kappa(A) = 10^5\) generated using the MATLAB commands \texttt{rng(1)} and \texttt{gallery('randsvd', [20, 20], 1e5, 1)}.
    The vector \(b\) is selected as the right singular vector corresponding to the fourth largest singular value, and the initial guess \(x_0\) is the zero vector.
    
    For this specific linear system, we use the restarted \(s\)-step GMRES with monomial\slash{}Newton\slash{}Chebyshev basis and BCGSI+.
    For simplicity, we do not consider the preconditioned version.
    We consider an extreme scenario where \(s\)-step GMRES undergoes a restart every \(20\) iterations, which corresponds to the dimension of \(A\).
    Then it can be ensured that the condition number of each sub-block \(\tilde B_{(i-1)s:is}\) is less than \(1.0 \cdot 10^5\) by setting \(s = 3\).
    However, \(\kappa(\tilde B_{1:p}) > 10^9\), and the relative backward error \(\frac{\norm{Ax - b}}{\normF{A} \norm{x} + \norm{b}}\) of the solution computed by the \(s\)-step GMRES with monomial\slash{}Newton\slash{}Chebyshev basis is at best around \(10^{-8}\), even when \(p = 20\), as shown in Figure~\ref{fig:badcase}.
    Note that for \(s = 3\), the backward error can be refined by the restart process.
    
    Unfortunately, the restart process does not always work, for example, \(s=4\) as shown in Figure~\ref{fig:badcase}.
    For this case, the condition number of each sub-block \(\tilde B_{(i-1)s:is}\) is less than \(1.0 \cdot 10^{10}\), but the relative backward error is at best around \(10^{-5}\) even using restart process.
    This means that we cannot control the condition number of the entire basis through controlling the condition number of each sub-block for the basis, and cannot use a restart process to recover the backward stability. 
\end{example}

Observe that \(\hat B_{1:p} = \hat V_{1:p}\) is evidently well-conditioned when \(s = 1\), and the issue described above arises only when \(s > 1\).
Therefore, to overcome this problem, we aim to ensure that \(B_{1:p}\) is near orthonormal, making \(B_{1:p}\) well-conditioned.
However, directly computing a QR factorization of \(K_{1:p}\) is not feasible because \(K_{1:p}\) becomes very ill-conditioned as \(s\) increases.
Note that \(\Span(V_{1:j}) = \Span (B_{1:j})\) as established in Theorem~\ref{thm:krylovsp}, and \(V_{1:j}\) remains well-conditioned until convergence.
Thus instead, we compute a QR factorization of the matrix \(\bmat{V_{1:(i-2)s+1:(i-1)s} & K_{(i-1)s+1:is}}\) during the \(i\)-th iteration, based on the orthonormal matrix \(V_{1:(i-2)s+1:(i-1)s}\).
This implies that \(B_{(i-1)s+1:is}\) is selected to be the \(Q\)-factor from the QR factorization of \((I - V_{1:(i-1)s} V_{1:(i-1)s}\trans) K_{(i-1)s+1:is}\).

The resulting \emph{modified} \(s\)-step Arnoldi process is outlined in Algorithm~\ref{alg:sstep-Arnoldi-modified}.
The distinction between the classical and modified \(s\)-step Arnoldi algorithms lies in the fact that the modified version employs an extra QR factorization, which can be computed by low-synchronization QR algorithms, to obtain \(B_{(i-1)s+1:is}\) in Line~\ref{line:mArnoldi:B} of Algorithm~\ref{alg:sstep-Arnoldi-modified}, rather than directly using \(K_{(i-1)s+1:is}\) as \(B_{(i-1)s+1:is}\) in Line~\ref{line:Arnoldi:B} of Algorithm~\ref{alg:sstep-Arnoldi-classical}.
This modified approach allows for the utilization of a significantly larger \(s\), with the trade-off of increasing the computation and communication cost of the QR factorization almost twofold; note that for sufficiently large $s$, this still provides an asymptotic communication savings versus standard GMRES.

\begin{algorithm}[!tb]
\begin{algorithmic}[1]
    \caption{The \(i\)-th step of the modified \(s\)-step Arnoldi process \label{alg:sstep-Arnoldi-modified}}
    \REQUIRE
    A matrix \(A \in \mathbb R^{n\times n}\), a vector \(r\), a block size \(s\), a left-preconditioner \(M_L \in \mathbb R^{n\times n}\), a right-preconditioner \(M_R \in \mathbb R^{n\times n}\),
    the basis \(B_{1:(i-1)s}\) and the preconditioned basis \(Z_{1:(i-1)s}\) generated by the first \(i-1\) classical Arnoldi steps,
    the matrix \(W_{1:(i-1)s}\), the orthonormal matrix \(V_{1:(i-1)s+1}\), and the upper triangular matrix \(R_{1:(i-1)s+1}\) satisfying \(\bmat{r & W_{1:(i-1)s}} = V_{1:(i-1)s+1} R_{1:(i-1)s+1}\).
    \ENSURE
    The basis \(B_{1:is}\), the preconditioned basis \(Z_{1:is}\), the matrices \(W_{1:is}\), \(V_{1:is+1}\), and \(R_{1:is+1}\) satisfying \(\bmat{r & W_{1:is}} = V_{1:is+1} R_{1:is+1}\).

    \STATE \(K_{(i-1)s+1:is} \gets \bmat{p_0(A) V_{(i-1)s+1} & p_1(A) V_{(i-1)s+1} & \dotsi & p_{s-1}(A) V_{(i-1)s+1}}\).
    \STATE \(B_{(i-1)s+1:is}\) is the \(Q\)-factor of \((I - V_{1:(i-1)s} V_{1:(i-1)s}\trans)^2 K_{(i-1)s+1:is}\) satisfying \((I - V_{1:(i-1)s} V_{1:(i-1)s}\trans)^2 K_{(i-1)s+1:is} = B_{(i-1)s+1:is} S_{(i-1)s+1:is}\). \label{line:mArnoldi:B}
    \STATE \(Z_{(i-1)s+1:is} \gets M_R^{-1} B_{(i-1)s+1:is}\).
    \STATE \(W_{(i-1)s+1:is} \gets M_L^{-1} A Z_{(i-1)s+1:is}\).
    \STATE Compute the QR factorization of \(\bmat{r & W_{1:is}} = V_{1:is+1} R_{1:is+1}\) based on \(\bmat{r & W_{1:(i-1)s}} = V_{1:(i-1)s+1} R_{1:(i-1)s+1}\). \label{line:modified-orth}
\end{algorithmic}
\end{algorithm}

As the above discussion, we first prove that the space spanned by \(B_{1:is}\) is the same as the space spanned by \(K_{1:is}\), the Krylov basis, in exact arithmetic.

\begin{theorem} \label{thm:krylovsp}
    Assume that \(K_{1:is}\) is defined by~\eqref{eq:def-Zi}, and \(B_{1:is}\), \(V_{1:is+1}\) are obtained via Algorithm~\ref{alg:sstep-GMRES} with Algorithm~\ref{alg:sstep-Arnoldi-modified}.
    Then
    \begin{align}
        & \Span\{B_{1:is}\} = \Span\{K_{1:is}\} = \Span\{V_{1:is}\} = \Span\bigl\{\bmat{r & A B_{1:is-1}}\bigr\}.\label{eq:thm:krylov:spanB=spanV}
    \end{align}
\end{theorem}

\begin{proof}
    We prove this theorem by induction.
    For the base case, \(B_{1:s}\) is the \(Q\)-factor of \(K_{1:s}\), i.e., \(K_{1:s} = B_{1:s} S_{1:s}\) and \(K_{1:s-1} = B_{1:s-1} S_{1:s-1}\), which amounts to
    \begin{align*}
        & \Span\{B_{1:s-1}\} = \Span\{K_{1:s-1}\},
        \qquad \Span\{B_{1:s}\} = \Span\{K_{1:s}\}.
    \end{align*}
    Together with the definition~\eqref{eq:def-Zi} of \(K_{1:s}\), we derive
    \begin{equation}
        \begin{split}
            \Span\{V_{1:s}\} ={}& \Span\bigl\{\bmat{r & A B_{1:s-1}}\bigr\} \\
            ={}& \Span\bigl\{\bmat{r & A K_{1:s-1}}\bigr\} \\
            ={}& \Span\bigl\{\bmat{r & A p_0(A) r & \dotsi & A p_{s-1}(A) r}\bigr\} \\
            ={}& \Span\{K_{1:s}\},
        \end{split}
    \end{equation}
    which also implies \(\Span\bigl\{\bmat{r & A B_{1:s-1}}\bigr\} = \Span\bigl\{\bmat{r & A K_{1:s-1}}\bigr\} = \Span\{V_{1:s}\}\).
    
    Then assuming that these hold for \(j-1\), i.e., 
    \begin{align}
        & \Span\{B_{1:(j-1)s}\} = \Span\{K_{1:(j-1)s}\} = \Span\{V_{1:(j-1)s}\}, \\
        & \Span\bigl\{\bmat{r & A B_{1:(j-1)s-1}}\bigr\} = \Span\bigl\{\bmat{r & A K_{1:(j-1)s-1}}\bigr\} = \Span\{V_{1:(j-1)s}\}, \label{eq:lem-krylovsp:spanZV}
    \end{align}
    we aim to prove that these hold for \(j\).
    Recalling Algorithm~\ref{alg:sstep-Arnoldi-modified}, \(B_{(j-1)s+1:js}\) is the \(Q\)-factor of \((I - V_{1:(j-1)s} V_{1:(j-1)s}\trans)^2 K_{(j-1)s+1:js}\) in exact arithmetic.
    Thus, we derive
    \begin{equation} \label{eq:lem-krylovsp:spanB}
        \begin{split}
            \Span\{B_{1:js}\}
            ={}& \Span\{\bmat{B_{1:(j-1)s} & B_{(j-1)s+1:js}}\} \\
            ={}& \Span\{\bmat{V_{1:(j-1)s} & (I - V_{1:(j-1)s} V_{1:(j-1)s}\trans)^2 K_{(j-1)s+1:js}}\} \\
            ={}& \Span\{\bmat{V_{1:(j-1)s} & K_{(j-1)s+1:js}}\} \\
            ={}& \Span\{K_{1:js}\}.
        \end{split}
    \end{equation}
    By \((I - V_{1:(j-1)s} V_{1:(j-1)s}\trans)^2 K_{(j-1)s+1:js-1} = B_{(j-1)s+1:js-1} S_{1:js-1, 1:js-1}\), we similarly have \(\Span\{B_{1:js-1}\} = \Span\{K_{1:js-1}\}\) and further from \(\bmat{r & A B_{1:js-1}} = V_{1:js} R_{1:js, 1:js}\)
    \begin{equation} \label{eq:lem-krylovsp:spanrAB-js-1}
        \Span\{V_{1:js}\} = \Span\bigl\{\bmat{r & A B_{1:js-1}}\bigr\} = \Span\bigl\{\bmat{r & A K_{1:js-1}}\bigr\}.
    \end{equation}
    Thus, by~\eqref{eq:lem-krylovsp:spanZV} and~\eqref{eq:lem-krylovsp:spanrAB-js-1}, we obtain
    \begin{equation}
        \begin{split}
            & \Span\{V_{1:js}\} \\
            &\quad = \Span\bigl\{\bmat{r & A K_{1:js-1}}\bigr\} \\
            &\quad = \Span\{V_{1:(j-1)s}, A K_{(j-1)s:js-1}\} \\
            &\quad = \Span\{V_{1:(j-1)s}, A K_{(j-1)s}, A p_0(A) V_{(j-1)s+1}, \dotsc, A p_{s-2}(A) V_{(j-1)s+1}\} \\
            &\quad = \Span\{V_{1:(j-1)s}, V_{(j-1)s+1}, A p_0(A) V_{(j-1)s+1}, \dotsc, A p_{s-2}(A) V_{(j-1)s+1}\} \\
            &\quad = \Span\{K_{1:(j-1)s}, K_{(j-1)s+1:js}\} \\
            &\quad = \Span\{K_{1:js}\}.
        \end{split}
    \end{equation}
    Thus, by induction on \(j\), we draw the conclusion~\eqref{eq:thm:krylov:spanB=spanV}.
\end{proof}

As described in Theorem~\ref{thm:b-Axi-H}, \(\kappa(\tilde B_{1:is})\), i.e., \(\kappa(\hat B_{1:is} D_{1:is})\) for any diagonal \(D_{1:is}\) with positive elements, has a critical influence on the backward error.
Note that we perform an extra QR factorization aiming to make \(B_{1:is}\) nearly orthonormal.
Thus, we consider \(D_{1:is}\) to be the identity matrix here.
In the following lemma, we give some inspiration to show why \(\kappa(\hat B_{1:is})\) can be expected to be well-conditioned regarding rounding errors.
The proof of the lemma can be found in Appendix~\ref{appendix:proof-lem-kappa-B}.

\begin{lemma} \label{lem:kappa-B}
    Assuming that \(V_{1:ks}\) is the exact result of Algorithm~\ref{alg:sstep-GMRES} with Algorithm~\ref{alg:sstep-Arnoldi-modified}, then there exists \(Y_{(i-1)s+1:is}\) such that
    \begin{equation} \label{eq:thm:krylov:B=VY}
        B_{(i-1)s+1:is} = V_{(i-1)s+1:is} Y_{(i-1)s+1:is}, \quad \forall i \leq k.
    \end{equation}
    
    Furthermore, assume that \(\hat B_{(i-1)s+1:is}\) and \(\hat V_{(i-1)s+1:is}\) are the computed results of Algorithm~\ref{alg:sstep-GMRES} with Algorithm~\ref{alg:sstep-Arnoldi-modified}.
    If there exists a small perturbation \(\Delta B_{(i-1)s+1:is}\) for any \(i \leq k\) such that
    \begin{equation} \label{eq:thm-kappaB:assump}
        \hat B_{(i-1)s+1:is} + \Delta B_{(i-1)s+1:is} = \hat V_{(i-1)s+1:is} \tilde Y_{(i-1)s+1:is},
    \end{equation}
    and it holds that
    \begin{equation} \label{eq:thm-kappaB:assump-2}
        \sum_{i=1}^k \omega_{B_i} + 5 \sqrt{ks}\, \sum_{i = 1}^k \normF{\Delta B_{(i-1)s+1:is}} + 7 s\, \omega_{k} \leq \frac{1}{2},
    \end{equation}
    then
    \(
        \kappa (\hat B_{1:ks}) \leq 2 \sqrt{n} + \sqrt{s},
    \)
    where \(\omega_{B_i}\) satisfies \(\normF{\hat B_{(i-1)s+1:is}\trans \hat B_{(i-1)s+1:is} - I} \leq \omega_{B_i}\).
\end{lemma}

\section{Numerical Experiments}
\label{sec:experiments}
In this section, we present numerical experiments to show that Algorithm~\ref{alg:sstep-GMRES} with Algorithm~\ref{alg:sstep-Arnoldi-modified} ($s$-step GMRES with the modified $s$-step Arnoldi process) can employ a much larger block size \(s\) compared to using the classical $s$-step Arnoldi process in Algorithm~\ref{alg:sstep-Arnoldi-classical}.
All tests are performed in MATLAB R2023a.

\subsection{Experiment settings}
The following variants of the \(s\)-step GMRES algorithms are tested:
\begin{enumerate}
    \item {\bf Classical \(s\)-step GMRES}: $s$-step GMRES (Algorithm~\ref{alg:sstep-GMRES}) with the classical $s$-step Arnoldi process (Algorithm~\ref{alg:sstep-Arnoldi-classical}), with stopping criterion~\eqref{eq:stopping-backward-err}.
    \item {\bf Modified \(s\)-step GMRES with the additional criterion}:  $s$-step GMRES (Algorithm~\ref{alg:sstep-GMRES}) with the modified Arnoldi process (Algorithm~\ref{alg:sstep-Arnoldi-modified}), with stopping criteria~\eqref{eq:stopping-backward-err} and~\eqref{eq:stoppingg-H}.
    \item {\bf Classical \(s\)-step GMRES with the additional criterion}: $s$-step GMRES (Algorithm~\ref{alg:sstep-GMRES}) with the classical Arnoldi process (Algorithm~\ref{alg:sstep-Arnoldi-classical}), with stopping criteria~\eqref{eq:stopping-backward-err} and~\eqref{eq:stoppingg-H}.
\end{enumerate}
In these three variants, BCGSI+ is employed as the orthogonalization method, respectively, in Line~\ref{line:classical-orth} of Algorithm~\ref{alg:sstep-Arnoldi-classical} and in Line~\ref{line:modified-orth} of Algorithm~\ref{alg:sstep-Arnoldi-modified}.
Since the condition number of the monomial basis grows exponentially with \(s\), the Newton and Chebyshev bases are used to generate each sub-matrix defined in~\eqref{eq:def-Zi}.
To clearly demonstrate the theoretical results, we do not consider restarting and preconditioners in our numerical experiments, but we reiterate that one could extend the theoretical results to restarted \(s\)-step GMRES.
The thresholds \texttt{tol} in~\eqref{eq:stopping-backward-err} and \texttt{tolH} in~\eqref{eq:stoppingg-H} are set to, respectively, \(n \macheps\) and \(\sqrt{n} \macheps\).

For constructing the linear systems \(A x = b\), we choose three commonly-used sparse square matrices, shown in Table~\ref{table:info_of_matrices}, from the SuiteSparse
Matrix Collection%
\footnote{\url{https://sparse.tamu.edu}}
and the Matrix Market%
\footnote{\url{https://math.nist.gov/MatrixMarket}}
as \(A\).
The two with larger condition numbers from the Matrix Market are suggested as test problems in~\cite{CLS2024}.
The other matrix from the SuiteSparse
Matrix Collection is a relatively well-conditioned matrix.
The right-hand vector \(b\) is set to be the vector of all ones and the initial guess \(x_0 = 0\).
We have selected these particular linear systems for demonstration purposes because the standard GMRES algorithm can converge relatively quickly even without the use of a preconditioner.

\begin{table}[!tb]
\centering
\caption{Properties of test matrices: the condition number in this table is estimated by the MATLAB command \texttt{svd}.}
\label{table:info_of_matrices}
\begin{tabular}{ccc}
\hline
Name& Size& Condition number \\
\hline
\texttt{494\_bus}      & \hphantom{0}494 & \(2.42\times 10^6\)\hphantom{0}\\
\texttt{fs1836}        & \hphantom{0}183 & \(1.74\times 10^{11}\)\\
\texttt{sherman2}      & 1,080 & \(9.64\times 10^{11}\)\\
\hline
\end{tabular}
\end{table}

\subsection{Tests for different block size \(s\)}
\begin{figure}[!tb]
    \centering
    \includegraphics[width=1.0\textwidth]{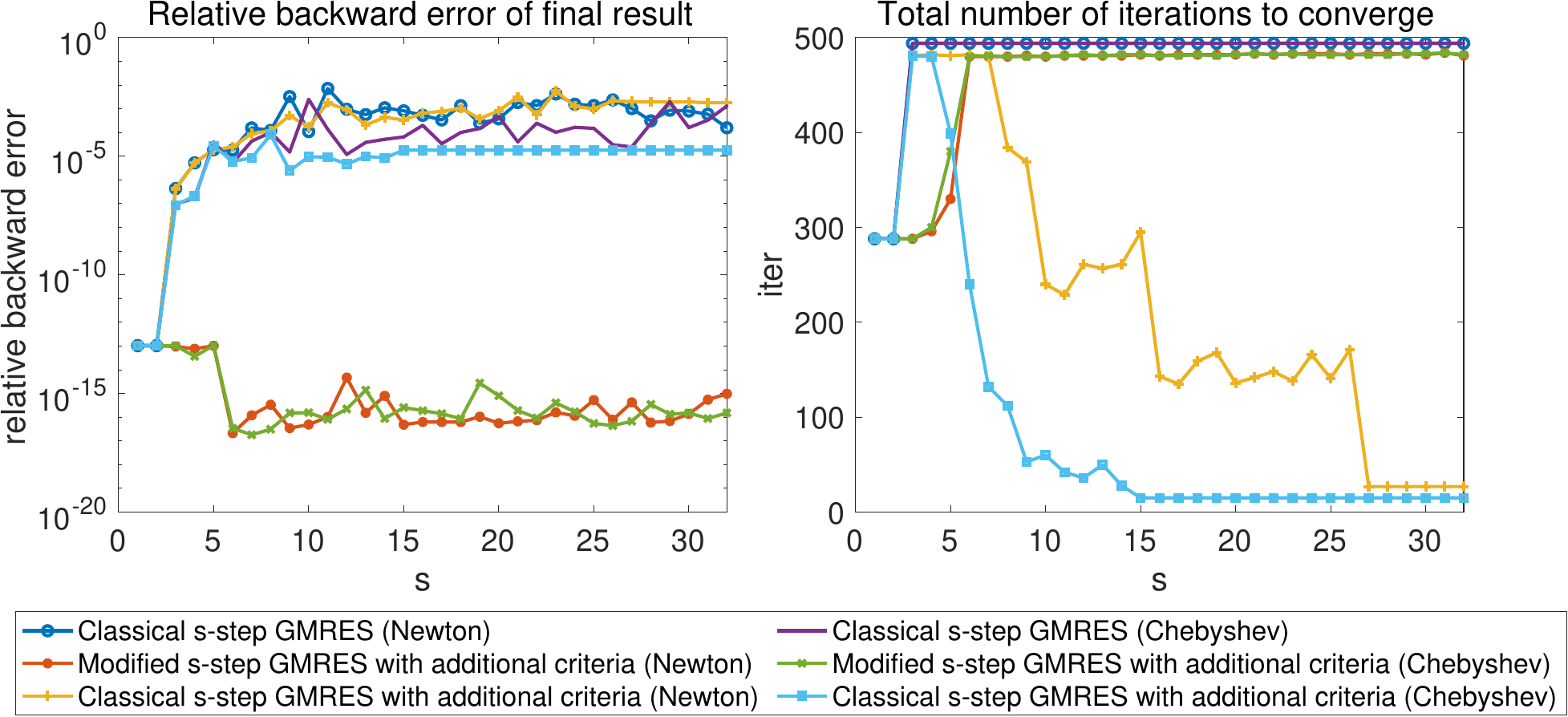}
    \caption{Relative backward errors (left) and the number of iterations (right) related to the block size \(s\), computed by different \(s\)-step GMRES algorithms, for \texttt{494bus}.}
    \label{fig:err-iter-s-494bus}
\end{figure}

\begin{figure}[!tb]
    \centering
    \includegraphics[width=1.0\textwidth]{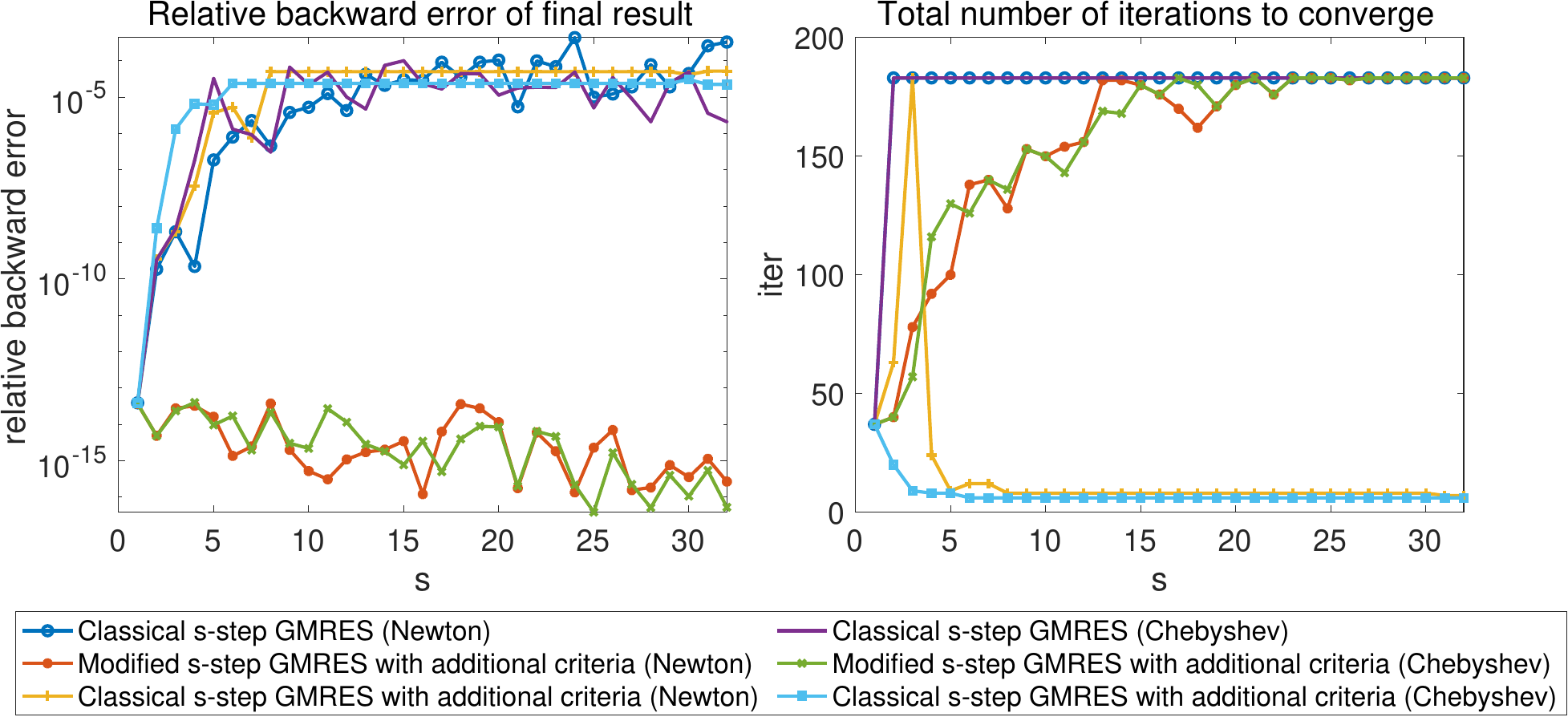}
    \caption{Relative backward errors (left) and the number of iterations (right) related to the block size \(s\), computed by different \(s\)-step GMRES algorithms, for \texttt{fs1836}.}
    \label{fig:err-iter-s-fs1836}
\end{figure}

\begin{figure}[!tb]
    \centering
    \includegraphics[width=1.0\textwidth]{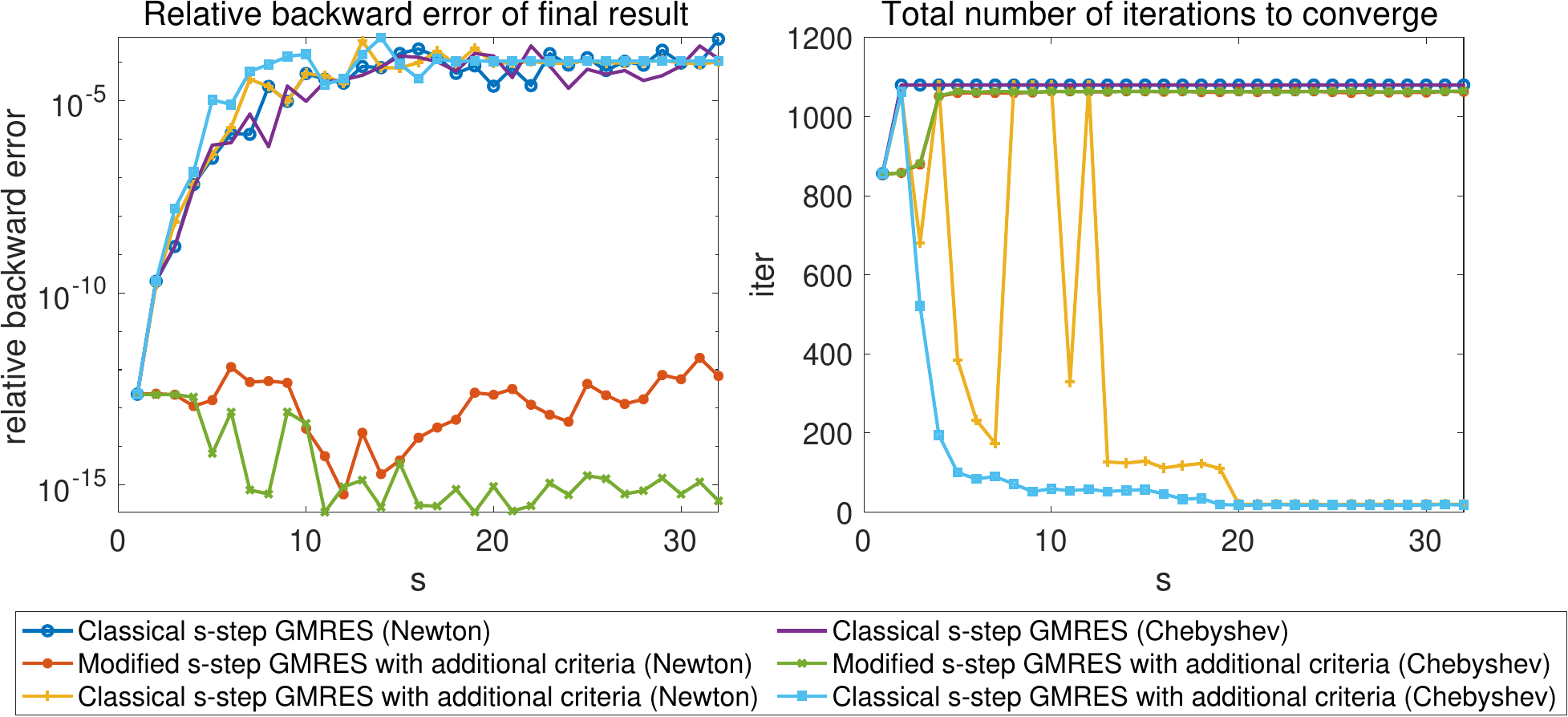}
    \caption{Relative backward errors (left) and the number of iterations (right) related to the block size \(s\), computed by different \(s\)-step GMRES algorithms, for \texttt{sherman2}.}
    \label{fig:err-iter-s-sherman2}
\end{figure}

We illustrate how varying the block size \(s\) impacts the relative backward error \(\frac{\norm{Ax - b}}{\normF{A} \norm{x} + \norm{b}}\) and the iteration count (in which each \(s\)-step counts as \(s\) iterations) for the different \(s\)-step GMRES variants in Figures~\ref{fig:err-iter-s-494bus}, \ref{fig:err-iter-s-fs1836}, and~\ref{fig:err-iter-s-sherman2}.

From these figures, it is clear that the classical \(s\)-step GMRES algorithm must use very small values of \(s\) to achieve satisfactory accuracy, since the condition number of \(B_{1:p}\) increases rapidly with larger \(s\).
It is further clear that the modified \(s\)-step GMRES algorithm benefits from the well-conditioned basis \(B_{1:p}\) and can utilize significantly larger \(s\) in practice without sacrificing accuracy.
It should be noted that until a point (depending on the matrix sparsity structure and machine parameters), a larger \(s\) results in lower communication cost, but it usually necessitates more iterations to reach a given level of backward error.
Therefore, it may not be beneficial in practice to set \(s\) too large, even if it appears that a larger \(s\) does not affect the backward error.

\begin{figure}[!tb]
    \centering
    \includegraphics[width=1.0\textwidth]{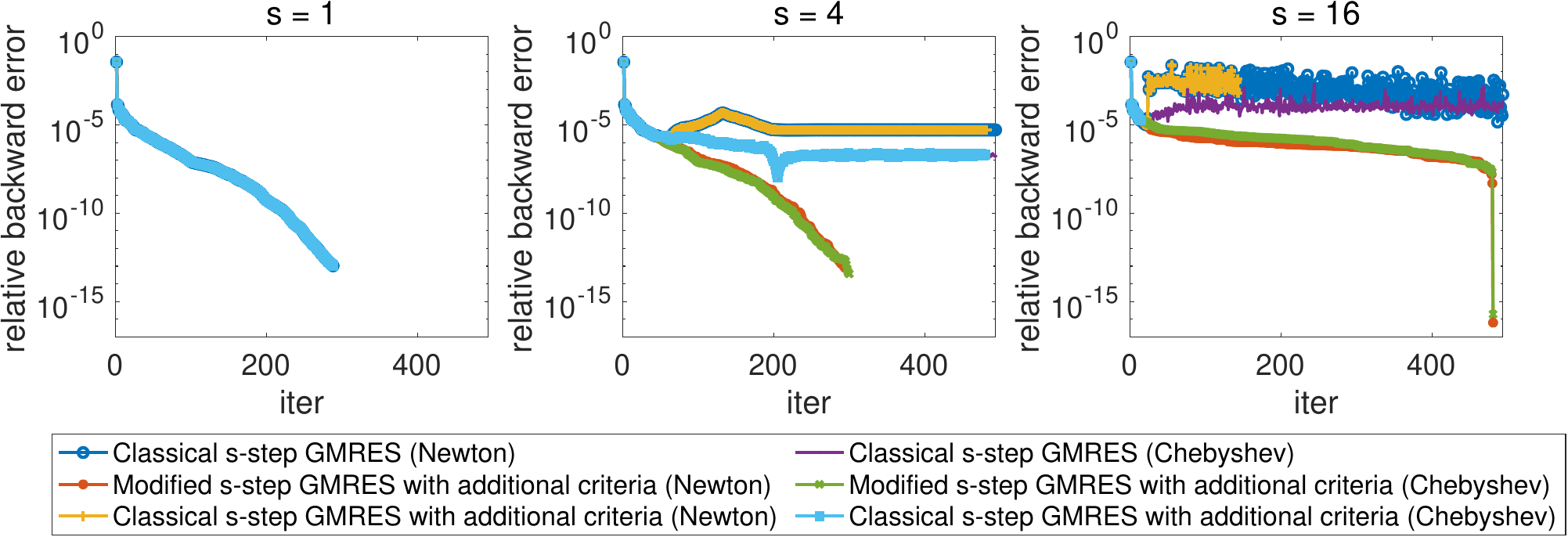}
    \caption{Relative backward errors by different \(s\)-step GMRES algorithms for \texttt{494bus} with \(s = 1\), \(4\), \(16\), respectively, from left to right.}
    \label{fig:err-iter-494bus}
\end{figure}

\begin{figure}[!tb]
    \centering
    \includegraphics[width=1.0\textwidth]{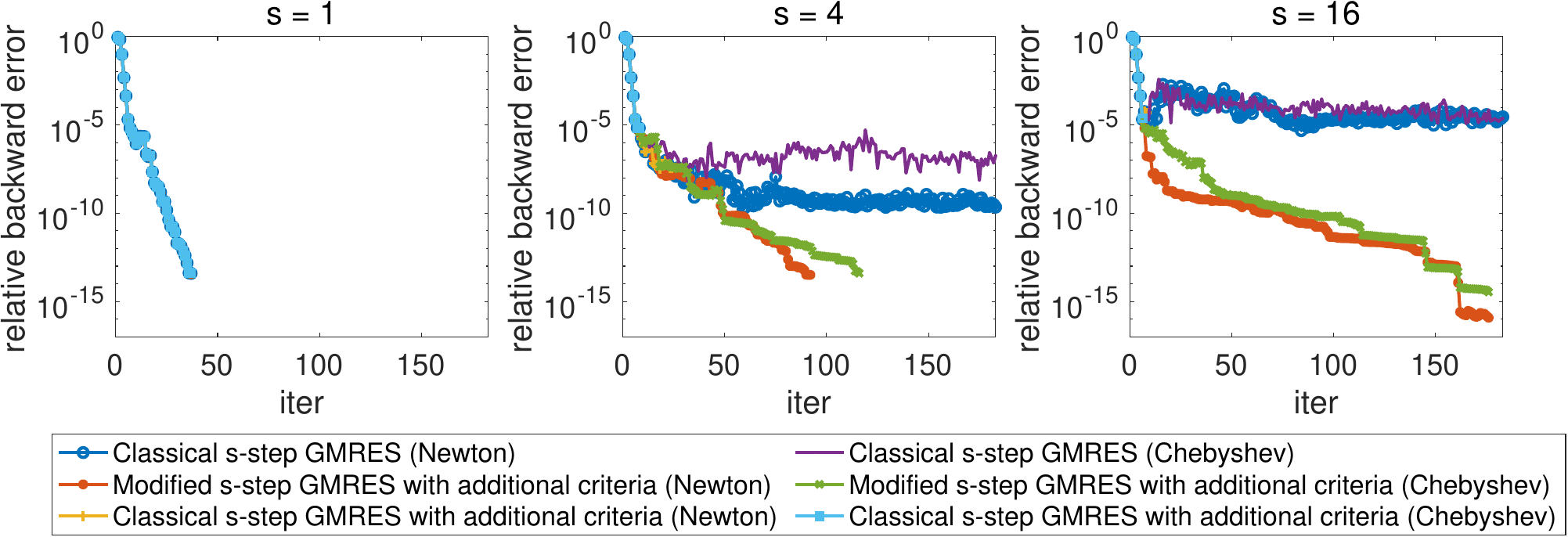}
    \caption{Relative backward errors by different \(s\)-step GMRES algorithms for \texttt{fs1836} with \(s = 1\), \(4\), \(16\), respectively, from left to right.}
    \label{fig:err-iter-fs1836}
\end{figure}

\begin{figure}[!tb]
    \centering
    \includegraphics[width=1.0\textwidth]{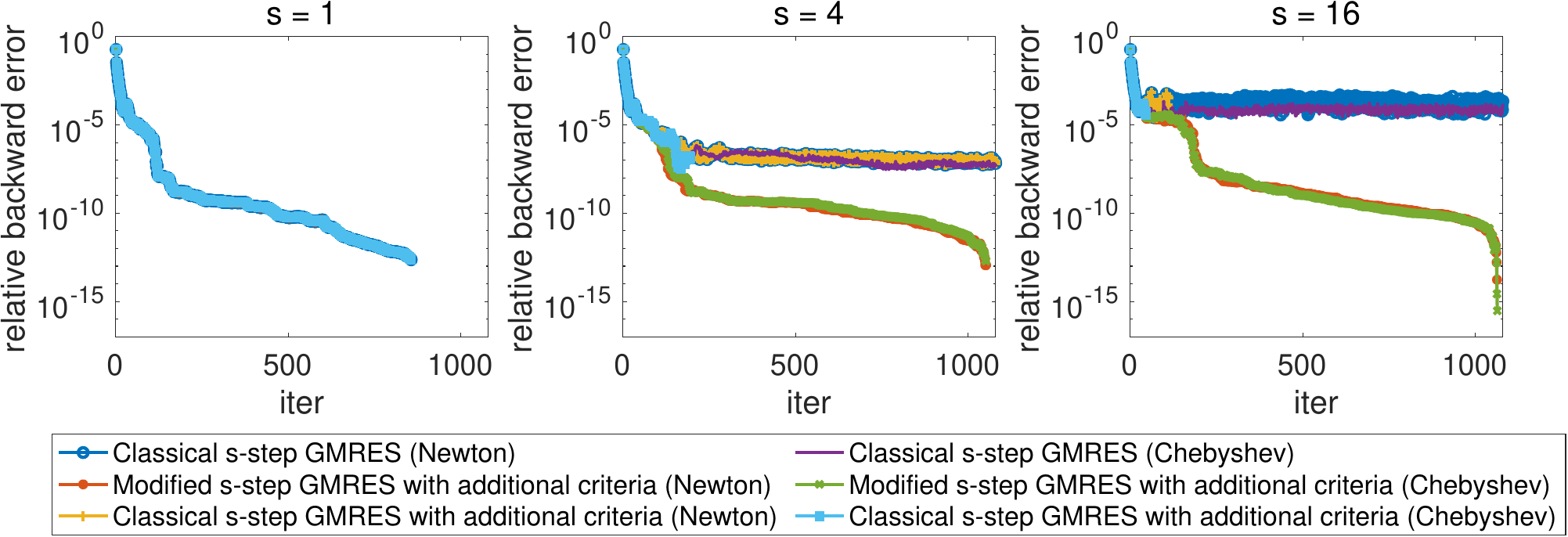}
    \caption{Relative backward errors by different \(s\)-step GMRES algorithms for \texttt{sherman2} with \(s = 1\), \(4\), \(16\), respectively, from left to right.}
    \label{fig:err-iter-sherman2}
\end{figure}

In Figures~\ref{fig:err-iter-s-494bus}, \ref{fig:err-iter-s-fs1836}, and~\ref{fig:err-iter-s-sherman2}, notice that the number of iterations suddenly decreases when using the additional criteria~\eqref{eq:stoppingg-H}.
To help clarify this observation,
we then choose three specific values of \(s\), i.e., \(s = 1\), \(4\), \(16\), to show the behavior of the backward error related to the iteration in Figures~\ref{fig:err-iter-494bus}, \ref{fig:err-iter-fs1836}, and~\ref{fig:err-iter-sherman2}.
As Lemma~\ref{lem:sstep-GMERS-bcgs2} predicts, there is almost no chance to obtain a better solution after the ``key dimension'' is reached; 
Figures~\ref{fig:err-iter-494bus}--\ref{fig:err-iter-sherman2} illustrate that using the criterion~\eqref{eq:stoppingg-H} achieves the "optimal" accuracy for most cases.
Thus, it is necessary to employ~\eqref{eq:stoppingg-H} to test the ``key dimension''.

\section{Conclusions}
\label{sec:conclusions}
In this work, we provide an abstract framework for analyzing the backward stability of the preconditioned \(s\)-step GMRES algorithm.
This framework accommodates various polynomial bases and block orthogonalization methods, demonstrating that the backward error of \(s\)-step GMRES is largely influenced by the condition number of the basis \(B\).
Moreover, it separately identifies the errors arising from the orthogonalization process (Step~2 in Section~\ref{sec:stability}) and from solving the least squares problem (Step~3), facilitating its application to \(s\)-step GMRES with different block orthogonalization techniques.
For standard GMRES, where \(s = 1\), our framework is an improved version of the modular GMRES framework introduced by~\cite{BHMV2024}.
We then apply the framework to analyze \(s\)-step GMRES with three widely-used block orthogonalization methods including block Householder QR, BMGS, and BCGSI+.
Furthermore, based on the framework, we give the stopping criteria and discuss the requirements of the orthogonalization methods used in GMRES, illustrating why the MGS, CGSI+, Householder QR, and TSQR algorithms are often used for orthogonalization in GMRES.

We then provide an example, specifically Example~\ref{example:badcase}, to illustrate the limitations of the classical \(s\)-step Arnoldi process (Algorithm~\ref{alg:sstep-Arnoldi-classical}).
In certain scenarios, the condition number of the basis \(B\) cannot be effectively managed by adaptively selecting \(s\) to control the condition number of each sub-block of the basis, utilizing Newton or Chebyshev bases, or a restart process.
Consequently, \(s\)-step GMRES may fail to achieve a satisfactory backward error even with a relatively small \(s\).
To address this issue, we introduce a modified \(s\)-step Arnoldi process that incorporates an extra QR factorization to make the basis \(B\) well-conditioned.  
Numerical experiments demonstrate that this modified \(s\)-step Arnoldi process allows for the use of a significantly larger \(s\) while achieving the required accuracy. While the modified approach increases the communication cost per iteration by a factor of 2, it may be beneficial in cases where a higher $s$ can lead to greater per-iteration speedup. Future work involves high performance implementations on large-scale problems in order to better evaluate these tradeoffs. 
\section*{Acknowledgments}
Both authors are supported by the European Union (ERC, inEXASCALE, 101075632). Views and opinions expressed are those of the authors only and do not necessarily reflect those of the European Union or the European Research Council. Neither the European Union nor the granting authority can be held responsible for them. The first author additionally acknowledges support from the Charles University Research Centre program No. UNCE/24/SCI/005.

\bibliographystyle{abbrvurl}
\bibliography{mybib}

\appendix
\section{Properties of BCGSI+}
\label{appendix:bcgsi+}
Given \(X \in \mathbb{R}^{m \times n}\) with \(m \geq n\), in Algorithm~\ref{alg:bcgs2}, we present the \((k+1)\)-th step of the BCGSI+ algorithm to compute \(X = QT\) with an orthonormal matrix \(Q \in \mathbb{R}^{m \times n}\) and an upper triangular matrix \(T \in \mathbb{R}^{n \times n}\).
Note that MGS or any unconditionally stable QR algorithm, e.g., Householder QR or Tall-Skinny QR (TSQR), described in~\cite{CLRT2022}, can be utilized in Line~\ref{line-bcgs2:qr1}, while any backward stable QR algorithm, i.e.,
\[
\hat W^{(2)} + \Delta W^{(2)} = \hat Q_{ks+1:(k+1)s} \hat T^{(2)}
\quad \text{with} \quad
\norm{\Delta W^{(2)}} \leq \bigO(\macheps) \norm{\hat W^{(2)}},
\]
can be employed in Line~\ref{line-bcgs2:qr2}; see~\cite{CLMO2024} for details.

\begin{algorithm}[!tb]
\begin{algorithmic}[1]
    \caption{The \((k+1)\)-th step of the BCGSI+ algorithm \label{alg:bcgs2}}
    \REQUIRE
    A matrix \(X \in \mathbb R^{m\times n}\), the block size \(s\), the orthogonal matrix \(Q_{1:ks} \in \mathbb{R}^{m \times ks}\) satisfying \(X_{1:ks} = Q_{1:ks} T_{1:ks, 1:ks}\).
    \ENSURE
    The orthogonal matrix \(Q_{1:(k+1)s} \in \mathbb{R}^{m \times (k+1)s}\) satisfies that \(X_{1:(k+1)s} = Q_{1:(k+1)s} T_{1:(k+1)s}\).

    \STATE \(S^{(1)} = Q_{1:ks}\trans X_{ks+1:(k+1)s} \in \mathbb{R}^{m \times s}\). \label{line:bcgs2:QTX}
    \STATE \(W^{(1)} = X_{ks+1:(k+1)s} - Q_{1:ks} S^{(1)}\). \label{line:bcgs2:W}
    \STATE Compute \(U \in \mathbb{R}^{m \times s}\) by QR algorithm such that \( W^{(1)} = U T^{(1)}\). \label{line-bcgs2:qr1}
    \STATE \(S^{(2)} = Q_{1:ks}\trans U \in \mathbb{R}^{m \times s}\).
    \STATE \(W^{(2)} = U - Q_{1:ks} S^{(2)}\).
    \STATE Compute \(Q_{ks+1:(k+1)s} \in \mathbb{R}^{m \times s}\) by QR algorithm such that \( W^{(2)} = Q_{ks+1:(k+1)s} T^{(2)}\). \label{line-bcgs2:qr2}
    \STATE \(T_{1:ks, ks+1:(k+1)s} = S^{(1)} + S^{(2)} T^{(1)}\).
    \STATE \(T_{ks+1:(k+1)s, ks+1:(k+1)s} = T^{(2)} T^{(1)}\).
\end{algorithmic}
\end{algorithm}

The properties of BCGSI+ have already been studied in~\cite[Section 2.3]{CLMO2024}.
Based on the results in~\cite[Section 2.3]{CLMO2024}, it is easy to obtain the following lemmas.

\begin{lemma} \label{lem:bcgs2:backwarderror-loo}
    Let \(\widehat Q_{1:js}\) and \(\widehat T_{1:js}\) be computed by Algorithm~\ref{alg:bcgs2}.
    If assuming \(\bigO(\macheps) \kappa(X_{1:js}) < 1\), then
    \begin{equation} \label{eq:bcgs2:backwarderror}
        X_{1:js} + \Delta X_{1:js} = \widehat Q_{1:js} \widehat T_{1:js}, \quad \norm{\Delta X_{i}} \leq \bigO(\macheps) \norm{X_i}
    \end{equation}
    for any \(i \leq js\), and
    \begin{equation} \label{eq:bcgs2:loo}
        \normF{\widehat Q_{1:js}\trans \widehat Q_{1:js} -  I} \leq \bigO(\macheps).
    \end{equation}
\end{lemma}

\begin{proof}
    Similarly to the proof of~\cite[Theorem 2]{CLMO2024}, it is easy to verify~\eqref{eq:bcgs2:backwarderror}, since each line of Algorithm~\ref{alg:bcgs2} is column-wise backward stable.
    The conclusion~\eqref{eq:bcgs2:loo} is directly followed by~\cite[Theorem 2]{CLMO2024}.
\end{proof}

From Lemma~\ref{lem:bcgs2:backwarderror-loo}, \(\epsqr = \bigO(\macheps)\), and~\eqref{eq:lem:LS:assump-rW} is satisfied when the \(Q\)-factor is not well-conditioned.
Then it remains to estimate \(\epsH\) defined by~\eqref{eq:lem-LS:assump-epsH}.

\begin{lemma} \label{lem:epsH}
    Assume that for \(is+j\), \(\widehat Q_{1:is+j}\) and \(\widehat T_{1:is+j}\) are computed by Algorithm~\ref{alg:bcgs2}.
    If
    \begin{align}
        & \normF{\widehat Q_{1:is+j-1}\trans \widehat Q_{1:is+j-1} - I} \leq \bigO(\macheps), \label{eq:lem-epsH:assump}\\
        & \normF{\widehat Q_{1:is+j}\trans \widehat Q_{1:is+j} - I} > \bigO(\macheps), \label{eq:lem-epsH:assump-2}
    \end{align}
    then
    \begin{equation}
        \max\{\norm{\widehat Q_{is+j} \widehat T_{is+j,is+j}}, \abs{\widehat T_{is+j,is+j}}\} \leq \bigO(\macheps) \norm{X_{is+j}}.
    \end{equation}
\end{lemma}

\begin{proof}
    Without loss of generality, we only need to prove the case \(j = 1\) since Algorithm~\ref{alg:bcgs2} is columnwise backward stable.
    From the assumption~\eqref{eq:lem-epsH:assump},
    \begin{equation}
        \normF{\widehat Q_{1:is+1}\trans \widehat Q_{1:is+1} - I}
        \leq \bigO(\macheps) + \normF{\widehat Q_{is+1}\trans \widehat Q_{is+1} - I} + 2\normF{\widehat Q_{1:is}\trans \widehat Q_{is+1}}.
    \end{equation}
    Note that~\cite[Equations (41)--(43), (46), (47), and (56)]{CLMO2024} do not depend on~\cite[Assumption (40)]{CLMO2024} when using an unconditionally stable QR algorithm in Line~\ref{line-bcgs2:qr1} of Algorithm~\ref{alg:bcgs2}.
    From~\cite[Equation (47)]{CLMO2024} it holds that
    \[
        \normF{\widehat Q_{is+1}\trans \widehat Q_{is+1} - I} \leq \bigO(\macheps).
    \]
    Then it remains to estimate \(\normF{\widehat Q_{1:is}\trans \widehat Q_{is+1}}\), which can be bounded as
    \begin{equation}
        \begin{split}
            \normF{\widehat Q_{1:is}\trans \widehat Q_{is+1}}
            & \leq \normF{\widehat Q_{1:is}\trans \tilde W^{(2)}_1 (\widehat T^{(2)}_{1,1})^{-1}} + \normF{\widehat Q_{1:is}\trans \Delta \tilde W^{(2)}_1 (\widehat T^{(2)}_{1,1})^{-1}} \\
            & \leq \frac{\normF{(I - \widehat Q_{1:is}\trans \widehat Q_{1:is}) \widehat Q_{1:is}\trans \widehat U_1}}{\abs{\widehat T^{(2)}_{1,1}}}
            + \frac{\normF{\widehat Q_{1:is}} \norm{\Delta \tilde W^{(2)}_1}}{\abs{\widehat T^{(2)}_{1,1}}},
        \end{split}
    \end{equation}
    where \(\widehat W^{(2)} = \tilde W^{(2)} + \Delta \tilde W^{(2)}\) with \(\tilde W^{(2)} = (I - \widehat Q_{1:is} \widehat Q_{1:is}\trans) \widehat U\).
    Together with~\cite[Equations (43) and (56)]{CLMO2024} and the assumption~\eqref{eq:lem-epsH:assump}, we obtain
    \begin{equation}
        \normF{\widehat Q_{1:is}\trans \widehat Q_{is+1}} \leq \frac{\bigO(\macheps)}{\abs{\widehat T^{(2)}_{1,1}}}.
    \end{equation}
    Then together with~\cite[Lemmas 2 and 6]{CLMO2024}, it follows that
    \[
        \tilde W^{(2)}_1 + \Delta W^{(2)}_1 = \widehat Q_{is+1} \widehat T^{(2)}_{1,1}, \quad
        \norm{\Delta W^{(2)}_1} \leq \bigO(\macheps),
    \]
    and further,
    \begin{equation*}
        \begin{split}
            \normF{\widehat Q_{1:is}\trans \widehat Q_{is+1}}
            & \leq \frac{\bigO(\macheps) \norm{\widehat Q_{is+1}}}{\norm{(I - \widehat Q_{1:is} \widehat Q_{1:is}\trans) \widehat U_1} - \norm{\Delta W^{(2)}_1}} \\
            & \leq \frac{\bigO(\macheps)}{\norm{(I - \widehat Q_{1:is} \widehat Q_{1:is}\trans) \widehat U_1} - \bigO(\macheps)}.
        \end{split}
    \end{equation*}
    This means that \(\normF{\widehat Q_{1:is+1}\trans \widehat Q_{1:is+1} - I} \leq \bigO(\macheps)\) if \(2 \norm{(I - \widehat Q_{1:is} \widehat Q_{1:is}\trans) \widehat U_1} > \bigO(\macheps)\) holds.
    Furthermore, the contrapositive is that \(2 \norm{(I - \widehat Q_{1:is} \widehat Q_{1:is}\trans) \widehat U_1} \leq \bigO(\macheps)\) if \(\normF{\widehat Q_{1:is+1}\trans \widehat Q_{1:is+1} - I} > \bigO(\macheps)\) guaranteed by the assumption~\eqref{eq:lem-epsH:assump-2}.
    Then by~\cite[Lemma 6]{CLMO2024}, and
    \[
        \widehat T_{is+1,is+1} = \widehat T^{(2)}_{1,1} \widehat T^{(1)}_{1,1} + \Delta T_{is+1,is+1} \quad \text{with} \quad 
        \abs{\Delta T_{is+1,is+1}} \leq \bigO(\macheps) \abs{\widehat T^{(2)}_{1,1}} \abs{\widehat T^{(1)}_{1,1}},
    \]
    we have
    \begin{equation*}
        \begin{split}
            \norm{\widehat Q_{is+1} \widehat T_{is+1,is+1}} \leq \norm{\widehat Q_{is+1} \widehat T^{(2)}_{1,1}} \abs{\widehat T^{(1)}_{1,1}} + \bigO(\macheps) \norm{\widehat Q_{is+1}} \abs{\widehat T^{(2)}_{1,1}} \abs{\widehat T^{(1)}_{1,1}}
            \leq \bigO(\macheps) \norm{X_{is+1}},
        \end{split}
    \end{equation*}
    which also implies \(\abs{\widehat T_{is+1,is+1}} \leq \bigO(\macheps) \norm{X_{is+1}}\) by noticing \(\norm{\widehat Q_{is+1}} \geq 1 - \bigO(\macheps)\).
\end{proof}

\section{Example for classical \(s\)-step GMRES with BCGSPIPI+}
\label{appendix:example-pip}
\begin{figure}[!tb]
    \centering
    \includegraphics[width=0.6\textwidth]{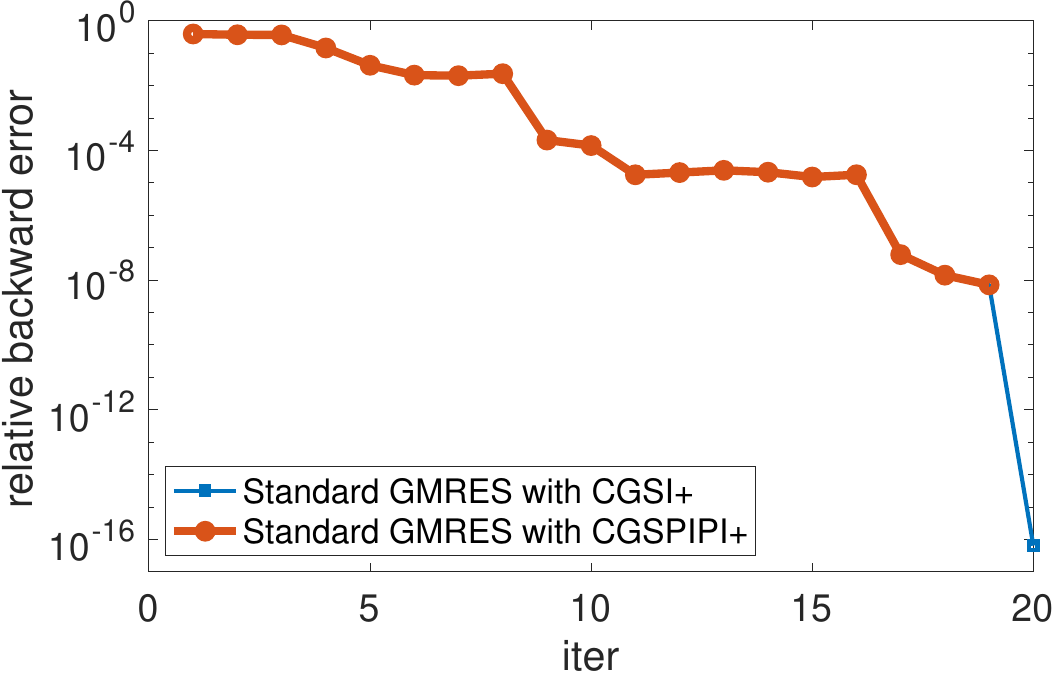}
    \caption{The plot for Example~\ref{example:pip}: The plot is of relative backward error. Each line with different color in the plot denotes the standard GMRES using different orthogonalization methods, including CGSI+ and CGSPIPI+.}
    \label{fig:pip-1}
\end{figure}

\begin{figure}[!tb]
    \centering
    \includegraphics[width=1.0\textwidth]{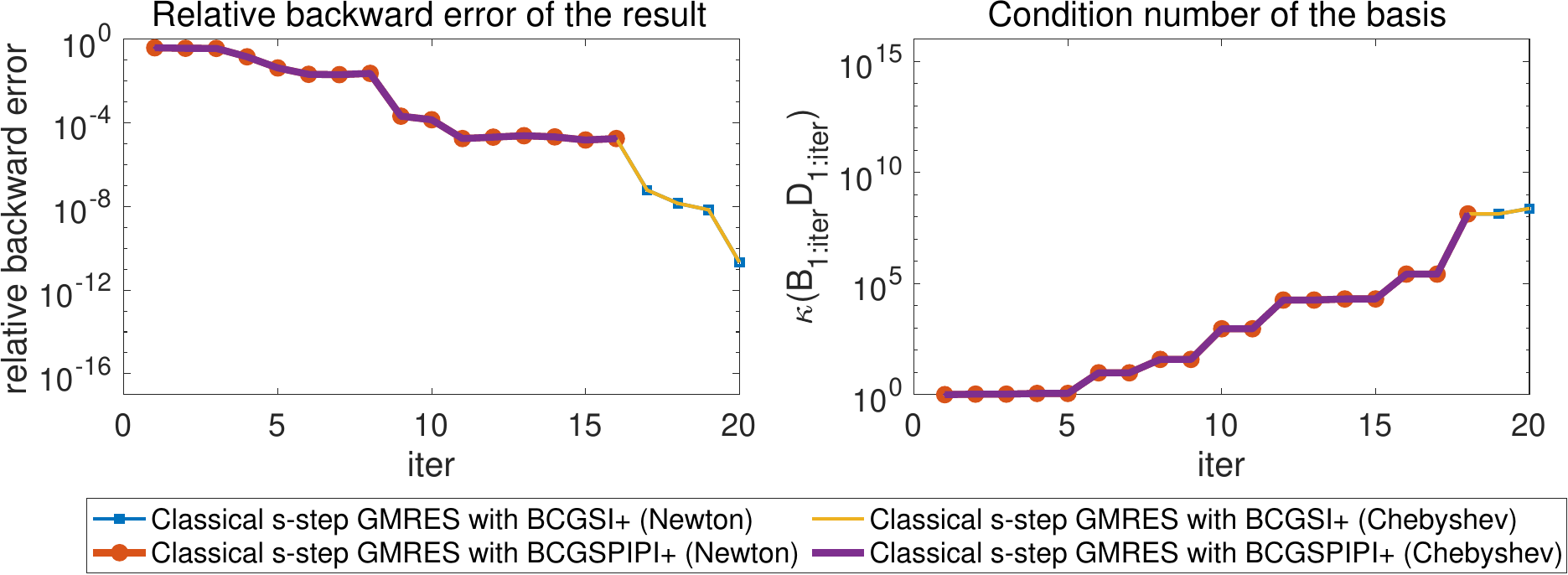}
    \caption{The plot for Example~\ref{example:pip}: From left to right, the plots are of relative backward error and the condition number of the basis \(\tilde{B}\), where we normalize each column of \(\hat B\) as \(\tilde B\).
    Each line with different color in the plots denotes the classical \(s\)-step GMRES with BCGSI+ or BCGSPIPI+ using different polynomials, including Newton and Chebyshev polynomials, to generate the basis by~\eqref{eq:def-Zi}.}
    \label{fig:pip}
\end{figure}

\begin{example} \label{example:pip}
    We construct the linear system \(Ax = b\), where \(A\) is a \(20\)-by-\(20\) random matrix with \(\kappa(A) = 10^{10}\) generated using the MATLAB command \texttt{rng(1)} and \texttt{gallery('randsvd', [20, 20], 1e10, 5)}.
    The vector \(b\) is selected as the right singular vector corresponding to the fourth largest singular value, and the initial guess \(x_0\) is the zero vector.
    
    For this specific linear system, the relative backward error \(\frac{\norm{Ax - b}}{\normF{A} \norm{x} + \norm{b}}\) of the solution computed by using standard GMRES with CGSPIPI+, namely BCGSPIPI+ with \(s = 1\), is approximately \(10^{-8}\).
    In contrast, using standard GMRES with CGSI+ results in an error of around \(10^{-16}\), as illustrated in Figure~\ref{fig:pip-1}.
    This difference occurs because CGSPIPI+ cannot generate a nearly orthonormal basis \(\hat V_{1:20}\) when the condition number of \(\bmat{\hat r & \hat W_{1:19}}\) exceeds approximately \(10^{8}\), implying \(\bigO(\macheps) \kappa^2(\bmat{\hat r & \hat W_{1:19}}) > 1\).
    Under this situation, the return value of \(\hat V_{20}\) is a \(\texttt{NaN}\) vector.
    Thus, it is not possible to obtain a more accurate solution than \(x^{(19)}\), whose backward error is approximately \(10^{-8}\).
    
    A similar result occurs for \(s\)-step GMRES with \(s = 2\).
    The relative backward error using BCGSPIPI+ is approximately \(10^{-5}\), while for \(s\)-step GMRES with BCGSI+, the error is around \(10^{-11}\), as illustrated in Figure~\ref{fig:pip}.
\end{example}

\section{Proof of Lemma~\ref{lem:kappa-B}}
\label{appendix:proof-lem-kappa-B}
\begin{proof}[Proof of Lemma~\ref{lem:kappa-B}]
    First, we aim to prove~\eqref{eq:thm:krylov:B=VY} by induction.
    For the base case, from~\eqref{eq:thm:krylov:spanB=spanV} with \(i=1\), there exists \(Y_{1:s}\) such that \(B_{1:s} = V_{1:s} Y_{1:s}\).
    Assume that \(B_{(i-1)s+1:is} = V_{(i-1)s+1:is} Y_{(i-1)s+1:is}\) holds for all \(i \leq j - 1\).
    Then our aim is to prove that it holds for \(j\).
    Since \(V_{1:(j-1)s}\) is orthonormal and \(B_{(j-1)s+1:js}\) is the \(Q\)-factor of \((I - V_{1:(j-1)s} V_{1:(j-1)s}\trans)^2 K_{(j-1)s+1:js}\), we obtain \(V_{1:(j-1)s}\trans B_{(j-1)s+1:js} = 0\).
    Together with the above assumptions on \(i \leq j - 1\) and~\eqref{eq:thm:krylov:spanB=spanV}, there exists \(Y_{(j-1)s+1:js}\) such that
    \begin{equation*}
        B_{1:js} = V_{1:js} \bmat{Y_{1:(j-1)s} & V_{1:(j-1)s}\trans B_{(j-1)s+1:js} \\ 0 & Y_{(j-1)s+1:js}}
        = V_{1:js} \bmat{Y_{1:(j-1)s} & 0 \\ 0 & Y_{(j-1)s+1:js}},
    \end{equation*}
    which gives~\eqref{eq:thm:krylov:B=VY} by induction on \(j\).

    Then we will bound \(\kappa (\hat B_{1:ks})\).
    By the definition of~\(\omega_{B_i}\), we only need to consider the off-diagonal blocks \(\normF{\hat B_{(i-1)s+1:is}\trans \hat B_{(j-1)s+1:js}}\).
    From~\eqref{eq:thm-kappaB:assump} and dropping the quadratic terms, it holds that
    \begin{equation*}
        \begin{split}
            & \normF{\hat B_{(i-1)s+1:is}\trans \hat B_{(j-1)s+1:js}} \\
            & \quad \leq \normF{\Delta B_{(i-1)s+1:is}} \norm{\hat V_{(j-1)s+1:js}} \norm{\tilde Y_{(j-1)s+1:js}} \\
            & \qquad + \norm{\tilde Y_{(i-1)s+1:is}} \normF{\hat V_{(i-1)s+1:is}\trans \hat V_{(j-1)s+1:js}} \norm{\tilde Y_{(j-1)s+1:js}} \\
            & \qquad + \normF{\Delta B_{(j-1)s+1:js}} \norm{\hat V_{(i-1)s+1:is}} \norm{\tilde Y_{(i-1)s+1:is}} \\
            & \quad \leq \normF{\Delta B_{(i-1)s+1:is}} \norm{\tilde Y_{(j-1)s+1:js}} \\
            & + \norm{\tilde Y_{(i-1)s+1:is}} \normF{\hat V_{(i-1)s+1:is}\trans \hat V_{(j-1)s+1:js}} \norm{\tilde Y_{(j-1)s+1:js}} \\
            & \qquad + \normF{\Delta B_{(j-1)s+1:js}} \norm{\tilde Y_{(i-1)s+1:is}},
        \end{split}
    \end{equation*}
    which implies that
    \begin{equation*}
        \begin{split}
            \sum_{i, j = 1; i \neq j}^k & \normF{\hat B_{(i-1)s+1:is}\trans \hat B_{(j-1)s+1:js}}^2 \\
            \leq{}& 3 \sum_{i, j = 1; i \neq j}^k \bigl(\normF{\Delta B_{(i-1)s+1:is}}^2 \norm{\tilde Y_{(j-1)s+1:js}}^2 \\
            & + \norm{\tilde Y_{(i-1)s+1:is}}^2 \normF{\hat V_{(i-1)s+1:is}\trans \hat V_{(j-1)s+1:js}}^2 \norm{\tilde Y_{(j-1)s+1:js}}^2 \\
            & + \normF{\Delta B_{(j-1)s+1:js}}^2 \norm{\tilde Y_{(i-1)s+1:is}}^2\bigr) \\
            \leq{}& 6k \sum_{i = 1}^k \bigl(\normF{\Delta B_{(i-1)s+1:is}}^2 \norm{\tilde Y_{(j-1)s+1:js}}^2\bigr) 
            + 3 \max_{i} \bigl(\norm{\tilde Y_{(j-1)s+1:js}}^4\bigr) \omega_k^2.
        \end{split}
    \end{equation*}
    Together with \(\normF{\hat B_{(i-1)s+1:is}} \leq \sqrt{s} (1 + \omega_{B_i})\) and
    \begin{equation}
        \norm{\tilde Y_{(i-1)s+1:is}} \leq \frac{\normF{\hat B_{(i-1)s+1:is}} + \normF{\Delta B_{(i-1)s+1:is}}}{1 - \normF{\hat V_{(i-1)s+1:is}\trans \hat V_{(i-1)s+1:is} - I}},
    \end{equation}
    we obtain
    \begin{equation*}
        \begin{split}
            & \normF{\hat B_{1:ks}\trans \hat B_{1:ks} - I}^2 \\
            & \quad = \sum_{i=1}^k \normF{\hat B_{(i-1)s+1:is}\trans \hat B_{(i-1)s+1:is} - I}^2 
            + \sum_{i, j = 1; i \neq j}^k \normF{\hat B_{(i-1)s+1:is}\trans \hat B_{(j-1)s+1:js}}^2 \\
            & \quad \leq \sum_{i=1}^k \omega_{B_i}^2
            + 6k \sum_{i = 1}^k \bigl(\normF{\Delta B_{(i-1)s+1:is}}^2 \norm{\tilde Y_{(j-1)s+1:js}}^2\bigr) + 3 \max_{i} \bigl(\norm{\tilde Y_{(j-1)s+1:js}}^4\bigr) \omega_k \\
            & \quad \leq \sum_{i=1}^k \omega_{B_i}^2
            + 24 ks \sum_{i = 1}^k \bigl(\normF{\Delta B_{(i-1)s+1:is}}^2\bigr) 
            + 48 s^2 \omega_k^2 
        \end{split}
    \end{equation*}
    by dropping the quadratic terms.
    This implies that, from the assumption~\eqref{eq:thm-kappaB:assump-2},
    \begin{equation}
        \begin{split}
            \normF{\hat B_{1:ks}\trans \hat B_{1:ks} - I} \leq \sum_{i=1}^k \omega_{B_i}
            + 5 \sqrt{ks} \sum_{i = 1}^k \normF{\Delta B_{(i-1)s+1:is}} + 7 s\, \omega_{k}
            \leq \frac{1}{2}.
        \end{split}
    \end{equation}
    Thus, \(\sigmin(\hat B_{1:ks})\) can be bounded by
    \begin{equation}
        \sigmin(\hat B_{1:ks}) \geq 1 - \normF{\hat B_{1:ks}\trans \hat B_{1:ks} - I} \geq \frac{1}{2}
    \end{equation}
    and further
    \begin{equation}
        \kappa (\hat B_{1:ks}) \leq \frac{\sum_i \normF{\hat B_{(i-1)s+1:is}}}{\sigmin(\hat B_{1:ks})}
        \leq 2 \sqrt{\sum_{i = 1}^k s\, (1 + \omega_{B_i})}
        \leq 2 \sqrt{n} + \sqrt{s}.
    \end{equation}
\end{proof}

\end{document}